\documentclass[12pt,leqno]{article}
\usepackage[latin1]{inputenc}
\usepackage{amsmath,amsfonts,amssymb,amsthm,amscd}
\date{}

\newcommand\Z{\mathbb{Z}}

\newcommand\A{\mathbb{A}}

\newcommand\C{\mathbb{C}}

\newcommand\N{\mathbb{N}}

\newcommand\R{\mathbb{R}}

\newcommand\Q{\mathbb{Q}}

\newcommand\fa{\mathfrak{a}}

\newcommand\pp{\mathfrak{p}}

\newcommand\Mm{\mathcal{M}}

\newcommand\Dd{\mathcal{D}}

\newcommand\Oo{\mathcal{O}}

\newcommand\CSs{\mathcal{S}}

\newcommand\Ff{\mathcal{F}}
\newcommand\Res{\mathrm{Res}}
\newcommand\ind{\mathrm{ind}}

\newcommand{\ba}{\backslash}
\newcommand\rg{\rightarrow}
\newcommand\lgr{\longrightarrow}   
\newcommand\Gal{\mathrm{Gal}}

\newcommand\Ind{\mathrm{Ind}}

\newcommand\hh{\hat{h}}

\newcommand{\gd}{\mathbf{d}}

\newcommand{\nequiv}{\equiv\kern-4.5mm/}

\def\adots{\mathinner{\mkern1mu\raise1pt\vbox{\kern7pt\hbox{.}}
\mkern2mu\raise4pt\hbox{.}
\mkern2mu\raise7pt\hbox{.}\mkern1mu}}

\newcommand{\intt}{\int\kern-4.5mm-}

\newcommand\Max{\mathrm{Max}}

\newtheorem{monlem}{{Lemma}}[section]
\newtheorem{madef}[monlem]{{Definition}}
\newtheorem{maprop}{{Proposition}}[section]
\newtheorem{montheo}{{Theorem}}[section]
\newtheorem{myhypo}{{Hypothesis}}[section]

\numberwithin{equation}{section}

\begin{document}

\title{A universal lower bound for certain quadratic integrals of automorphic $L$--functions}

\author{Laurent Clozel}

\maketitle

{\large With an Appendix by Laurent Clozel and Peter Sarnak}

\vskip5mm

 \textbf{\Large{Introduction}}
 
 \vskip3mm

\noindent\textbf{1.} Let $\pi$ be a unitary cuspidal representation of $GL(m,\A)$ where $\A$ is the ring of ad\`eles of $\Q$. Let $L(s,\pi)$ denote its $L$--function. More generally, we may consider a representation $\pi$ of the form $\pi_1\times \pi_2 \times \cdots \times \pi_r$ where $m=m_1+m_2+\cdots +m_r$, $\pi_i$ is a cuspidal unitary representation of $GL(m_i,\A)$ and $\times$ denotes induction by blocks. Thus
$$
L(s,\pi) = L(s,\pi_1)L(s,\pi_2) \cdots L(s,\pi_r).
$$

We will assume that $L(s,\pi)$ has at most a simple pole, at $s=1$. Thus there is at most one occurence of $m_i=1$, $\pi_i(x)=|x|^c$ $(x\in\A^\times)$ and we then assume that $c=0$, so $L(s,\pi_i)=\zeta(s)$. (We can in fact consider more general $L$--functions, coming from a finite extension $F/\Q$. See~\S~2.1)

Assume $s$ is a zero of $L(s,\pi)$ in the critical strip : $s=\sigma+i\tau$, $0<\sigma<1$. We first assume that $\pi$ contains no ``z\^eta'' factor. In this Introduction we assume that $s$ is on the critical line : $\sigma=1/2$.

\vspace{2mm}

\textbf{Theorem A.} ($\pi_i=\C$ \textit{does not occur} in $\pi$) $(\sigma=1/2)$
$$
\int_{-\infty}^{+\infty} \Big|\frac{L(\frac{1}{2}+it,\pi)}{\frac{1}{2}+it-s}\Big|^2 dt > 2\pi \log 2.
$$

In the case where  $L(s,\pi)$ has a z\^eta factor, the result is slightly  \\ different :

\vspace{2mm}

\textbf{Theorem B.} $(\sigma=1/2)$
$$
(i) \\ \int_{-\infty}^{+\infty} \Big|\frac{L(\frac{1}{2}+it,\pi)}{\frac{1}{2}+it-s}\Big|^2 dt > 2\pi \Big(\log 2 - \frac{2\vert \kappa \vert}{|1-s|}\Big)
$$
\textit{where $\kappa$ is the residue of $L(s,\pi)$ at $s=1$.}
$$
(ii) \\   \int_{-\infty}^{+\infty} \Big|\frac{L(\frac{1}{2}+it,\pi)}{\frac{1}{2}+it-s}\Big|^2 dt > \pi \log 2.
$$

\vspace{2mm}

Obviously Theorem B contains Theorem A. Note that $\kappa\ll D^\varepsilon$ for small $\varepsilon>0$, $D$ being the conductor of $\pi$ and the implicit constant being uniform (for $m$ fixed) if $\pi$ satisfies the Ramanujan hypothesis. See Iwaniec--Kowalski \cite[p.~160]{IK}. Thus the first lower bound is effective.

\vspace{2mm}

The two next results were suggested by Peter Sarnak. They correspond to $(s=0)$ in the previous statements; note that $L(s, \pi)$ does not necessarily vanish at $(s=0).$

\vspace{2mm}

\textbf{Theorem C.} ($\pi_i=\C$ \textit{does not occur} in $\pi$)
$$
\int_{-\infty}^{+\infty} \Big|\frac{L(\frac{1}{2}+it,\pi)}{\frac{1}{2}+it}\Big|^2 dt > \pi .
$$

\vskip2mm

This theorem implies a uniform  'pseudo- $\Omega$-result':

\vskip2mm

\noindent \textbf{{Corollary.}} \textit{Under the same assumptions, for any} $\varepsilon<1/2$,
$$ 
Sup \\_{t\in \R}  \frac{\vert L(\frac{1}{2}+it,\pi) \vert}{\vert \frac{1}{2}+it \vert^{\varepsilon}}  > \sqrt2/2 (1+ O(\epsilon))
$$
\textit{where the remainder in $O(\epsilon)$  is independent of $m$ and $\pi$.} 

\vskip2mm

Of course this is \textit{not} an $\Omega$-result, as we do not obtain the inequality for arbitrarily large values of $t$. Note that the Sup is finite according to the Lindel\"{o}f conjecture. The proof follows directly by bounding the integral by the product of the supremum (squared) and of the integral of $\vert \frac{1}{2}+it \vert^{-2+2\varepsilon}$. By considering the powers $L(s,\pi)^r$, $r \geq 1$ (cf. the argument in the next paragraph), we can even deduce:

\vskip2mm

\noindent \textbf{{Corollary.}} \textit{Under the same assumptions,} 
$$ 
lim _{\varepsilon \rightarrow0} \\ Sup \\_{t\in \R}  \frac{\vert L(\frac{1}{2}+it,\pi) \vert}{\vert \frac{1}{2}+it \vert^{\varepsilon}} \geq 1
$$
\textit{uniformy with respect to  $m$ and $\pi$.}

\vspace{2mm}

In the case where a residue is present, we obtain a slightly weaker result:

\vskip2mm

\textbf{Theorem D.} 
$$
\int_{-\infty}^{+\infty} \Big|\frac{L(\frac{1}{2}+it,\pi)}{\frac{1}{2}+it}\Big|^2 dt > \pi /2.
$$

\noindent\textbf{2.} In the present state of our knowledge, it is not known that the integrals in Theorems~A-D  are finite. One may consider the Theorems true if they are infinite: so construed they will be proved in the Appendix by Sarnak and this author. In view of classical results for $m=1$, and of recent subconvexity results for $m=2$, they are finite in these cases. \textit{The main body of the paper is concerned with the cases where all the integrals considered are finite.} We prove the Theorems for $m=1$ (Riemann z\^eta function  or Dirichlet $L$--series), see Theorem~1.2, Proposition~3.1\footnote{The constants in Proposition 3.1 are different because we do not assume $\sigma=1/2$. For $\sigma=1/2$ one gets the indicated constants.}. For $m=2$, this includes the case where $L(s,\pi)=\zeta_F(s)$ for $F$ a quadratic extension of $\Q$ ; and of $L(s,\pi)$ for a cuspidal (unitary) representation of $GL(2,\A_\Q)$. \textit{We assume} that the Archimedean factor $\pi_\infty$ is self--dual and that $\pi$ verifies the Ramanujan Conjecture. In this case Theorems A and B are proved in Proposition~3.2\footnote{Likewise, the constants can be replaced by the constant $2\pi\log 2$ of Conjecture~A if $Re (s)=1/2$.};  Theorems C and D are proved in \S 3.3.

Sarnak also pointed out that lower bounds, uniform in $m$, such as in Theorem C cannot be obtained for small intervals. Let $L(s,E)$ be the $L$-function of an elliptic curve $E$ over $\Q$, normalised as here: functional equation relating $L(E,s)$ and $L(E, 1-s)$. One can find $E$ such that $\vert L(\frac{1}{2}+it,E) \vert<1$ for $t\in[-1,1]$. See \cite{LF}, and in particular the $L$-function ~~ of ~~$X_0(11)$ there.\footnote{I thank Andrew Booker for this reference.}  Consider $L(s)=L(s,E)^r$. Then for $r$ sufficiently large the integral in Theorem C, restricted to$[-1, 1]$, is arbitrarily small. Note that the arguments in this paper do not require $\pi$ to be cuspidal.

\vspace{2mm}

\noindent\textbf{3.} The proof relies on an application, apparently new, of the Mellin transform. Consider
$$
L(s,\pi) = \sum_1^\infty a_n n^{-s}
$$
and, for $X>0$,
\setcounter{section}{0}
\begin{eqnarray}
A_s(X) &=& \displaystyle\sum_{n\le X} a_n n^{-s}\\
H_s(X)&=& X^{s-1} \displaystyle\sum_{n\le X} a_n n^{-s} - \frac{\kappa}{1-s}.
\end{eqnarray}
We were fist led (by Tate's thesis) to consider $H_s(X)$ when $\pi$ corresponds to the z\^eta function of a number field $F$ of degree $m$ over $\Q$ : in this case,
$$
H_s(X) =X^{s-1} \sum_{N\fa\le X} N\fa^{-s} - \frac{\kappa}{1-s}
$$
where $\fa$ ranges over integral ideals of $F$, different from $\{0\}$, and $\kappa$ is the usual residue. Note that if we consider $x \in F_\infty = \prod\limits_{v|\infty}F_v$ and $X=|x| = \prod\limits_{v|\infty} |x_v|$, $H_s(X)$ is a function of slow growth on $F_\infty$. In particular its Fourier transform  $\Ff H_s$ is defined.

In this case, $H_s$ enjoys some remarkable properties. Let
$$
K_s(X) = D^{-1/2}\ X^{s-1} \sum_{\fa \subset \Dd^{-1}\atop N\fa\le X} N\fa^{-s}-\kappa\frac{D^{1/2}}{1-s}.
$$
Here $\Dd^{-1}$ is the inverse different of the ring of integers of $F$ and $\fa$ ranges over fractional ideals, different from $0$. We now have\begin{equation}
\textit{For}\ 0<\sigma<1,\ \zeta_F(s)=0\ \textit{if,\ and\ only\ if,}\ \Ff(H_s) = - K_{1-s}.
\end{equation}

See Theorem 1.1.

Furthermore, using Perron's formula, one can compute the Mellin transform of $H_s$ :
$$
\Mm H_s(w)= \int_0^\infty H_s(x) x^{w-1}dx.
$$
One needs  an estimate for $H_s(x)$ $(x\rg \infty)$ ; as a first step we use an estimate of Landau (1915) for $A_0(x)$. Assuming again $\sigma=\frac{1}{2}$, one finds that $\Mm H_s$ is absolutely convergent for $0<Re(w)<\frac{1}{2}$ (Lemma~1.1). On the other hand, if $0<c<1$, $c\le\frac{2}{m}$, and $\zeta_F(s)=0$ :
$$
H_s^* (x)= \frac{1}{2i\pi} \int_{c-i\infty}^{c+i\infty} \frac{\zeta_F(1-w)}{1-s-w} x^{-w}dw.
$$
Here, as usual, $H_s^*(x)$ is given by $x^{s-1}\sum\limits_{n\le x} a_n n^{-s}$, the last term of the sum being weighted  by $\frac{1}{2}$ if $x=n$ is an integer.

This suggests, of course, that more generally :

\vspace{2mm}

\textbf{Conjecture E.}
$$
\Mm H_s (w) = \frac{L(1-w,\pi)}{1-s-w}
$$
\textit{when $H_s$ is defined by $(0.2)$, $L(s,\pi)=0$, and $\Mm H_s$ is defined.}

\vspace{2mm}

Returning to the case of $\zeta_F$, one finds that for $F$ quadratic $\frac{\zeta_F(1-w)}{1-s-w}$ is an $L^2$ function of $t$ $(w=1/2+it)$ and one can deduce (although $Re(w)=1/2$ is the limit of the domain of convergence of $\Mm H_s$) that
$$
\int_0^\infty |H_s(x)|^2 dx = \frac{1}{2\pi} \int_{1/2-i\infty}^{1/2+i\infty} \frac{|\zeta_F(w)|^2}{|w-s|^2}dw.
$$
This implies Theorem B by an easy computation of $\int_1^2 |H_s(x)|^2 dx$ (Theorem~1.2). The proof is easier for $\zeta=\zeta_\Q$ ; in this case $\Mm H_s$ converges for $0<Re\,w<1$ (Proposition~1.1, Theorem~1.2).

\vspace{2mm}

\noindent\textbf{4.} We now want to make sense, for the $L$--functions of \S~0.1, of Conjecture~E for $w$ in some domain which ideally should include the line $Re(w)=\frac{1}{2}$. For this we need a growth estimate for $H_s(x)$ when $x\rg +\infty$ and $L(s,\pi)=0$~; equivalently, for $A_s(x) -\frac{\kappa}{1-s} x^{1-s}$.

 As we pointed out, the remainder    $A_0(x) - \kappa x$ has been majorised by Landau, at least if $\pi$ verifies some conditions which are satisfied when $\pi$ is associated to $\zeta_F$ (see \S~2.3). For general $s$, we use the proof of a similar result by Friedlander and Iwaniec \cite{FI}. (They prove a stronger result, giving an expression of the remainder with an explicit constant depending on the conductor of~$\pi$.)

Friedlander and Iwaniec implicitly impose a condition on $\pi$, namely, that the Archimedean factor $\pi_\infty$ be self--dual ; furthermore, they assume that $\pi$ satisfies the Ramanujan Conjecture. (See \S~2.3.) We now make these assumptions. It is then possible to adapt their proof (for $A_0$) to the case of $A_s$ for $s$ in the critical strip. The final result is Theorem~2.2 and its Corollary : if $s$ is a zero,
$$
H_s(x) = O\Big(D^{\frac{1}{m+1}} x^{-\frac{2}{m+1}+\varepsilon}\Big).
$$
In particular, $H_s(x)$ is $L^2$ for $m\le2.$~\footnote{We have given rather abundant details in this part of the paper, as the article of Friedlander--Iwaniec is elliptic.}

Having done this, we can return to the conjecture for $m\le 2$. This is done in \S~3 : in this case Conjecture~E is well--defined, the functions on both sides are $L^2$, and we prove Theorem~B. Furthermore the vertical integrals of $\big|\frac{L(w,\pi)}{w-s}\big|^2$ can be majorized for some values of $Re(w)$ different for $\frac{1}{2}$ : see~\S~3.2. \S 3.3 is devoted to the case of $s=0$ (Theorems C,D.)

In \S 3.4, we pause to draw some consequences of the estimate of \S~2 for the abscissa of convergence of $L(s,\pi)$ (Theorem~3.1). The general result is that it converges for $Re(s) >  1-\frac{2}{m+1}$. We then discuss, in light of known conjectures,  the expected value of this abscissa.

The relevant Conjecture here is due to Friedlander--Iwaniec in the same paper \cite[Conjecture~2]{FI}. Assuming, as seems  likely, that their conjectural estimate can be extended, as in \S~2, from $(s=0)$ to an arbitrary value in the critical strip, it would imply (for $s$ a zero) that $H_s(x)\ll x^{-1/2-\frac{1}{2m}+\varepsilon}$. On the other hand, the integral in Theorem~B is, of course, finite if we assume an approximation of the Lindel\"of hypothesis. Under these assumptions, Conjecture~E would be meaningful and true, and the finite case of Theorem~B would follow. See \S~3.5.

One interesting aspect of the duality introduced by Conjecture~E is that it gives a relation between the growth of $L(s,\pi)$ in vertical lines (controlled by the Lindel\"of hypothesis) and the growth of $H_s$ (for which a useful control is given only by the Friedlander--Iwaniec Conjecture : the generalised Riemann hypothesis does not suffice.) See~3.5.

Finally, in \S 3.6, we hint at the following problem. We consider $\zeta=\zeta_\Q$ for simplicity ; in this case Theorem~B is true, the integral being finite : there is a lower bound for
$$
I(s) = \int_{-\infty}^{+\infty} \Big|\frac{\zeta(1/2+it)}{\frac{1}{2}+it-s}\Big|^2dt
$$
when $s$ is a zero. Does this integral tend to infinity with $s$ ? This seems likely ; a counter--intuitive consequence of the existence of a uniform bound for $I(s)$ is given in Proposition~3.3.

To conclude, we note that it would be desirable to rid the proof in \S~2 of the assumptions on $\pi$, but this seems difficult.  Moreover, let $\theta:0\le\theta\le1/2$ be the ``deviation from the Ramanujan hypothesis'' for $\pi$. (see e.g. \cite[\S~3.1]{Cl3}). Thus $\theta\le 1/2-\frac{1}{m^2+1}$ \cite{LRS} ; $\theta\le \frac{7}{32}$  for $m=2$ and Maass forms \cite{KS}. The (known) abscissa of convergence for $L(s,\pi)$ should be shifted  (positively) by $\theta$ ; similarly one can expect that $\Mm H_s$ is defined for $Re(w) <\frac{1}{2}+\frac{1}{2m}-\theta$. For Maass forms, this yields $Re(w)<\frac{1}{2}+\frac{1}{32}$. If so, one could consider the integral on $Re\,w=\frac{1}{2}$ and Theorem~A would be accessible with the same proof.

\vspace{2mm}

\textbf{Acknowledgements} :  I am very grateful to Peter Sarnak, for pointing out the unconditional form of the theorems and the case of $(s=0)$, and for suggesting that we write the Appendix. I thank R\'{e}gis de la Bret\`{e}che for helping me to read the paper of Friedlander-Iwaniec \cite{FI}. I also thank Andrew Booker, Emmanuel Kowalski, Herv\'{e} Jacquet, Etienne Fouvry, John Friedlander and J.-L. Waldspurger for useful correspondence.

\eject

\section{The Tate kernel for number fields, and its Fourier and Mellin transforms}

\subsection{} Let $F$ be a number field. We denote by $v$ a place of $F$. Let $F_\infty=\prod\limits_{v|\infty}F_v$. The Fourier transform on $F_\infty$ is defined place by place : on $F_v\cong \R$,
$$
\Ff f(y)=\hat{f}(y) = \int f(x){e}^{-2i\pi xy}dy.
$$

On $F_v\cong \C$, variables $z=x+iy$, $w=\xi+i\eta$,
$$
\Ff f(w) = \hat{f}(w) = \int f(z)e^{-4i\pi Re(zw)}dz
$$
where $dz =2dxdy$. (This is Tate's  self--dual normalisation, cf \cite{Ta}). This also defines the Fourier transform for distributions on $F_\infty$.

For $x=(x_v)\in F_\infty$, we write $X=\prod\limits_v |x_v|$ where the absolute value is normalised \cite{Ta}. Let $\Oo=\Oo_F$ be the integers of $F$, $\Dd=\Dd_F$ its different and $\Dd^{-1}$ the inverse different. Let $D=|D_F|$ the absolute value of the discriminant ; thus $D= |N_{F/\Q}\Dd|$.

We denote by $\fa$ a non--zero fractional ideal in $F$, and by $N\fa$ its norm. Let $\kappa=\frac{2^{r_1}(2\pi)^{r_2}hR}{w\sqrt{|D|}}$ be the residue at 1 of the z\^{e}ta function $\zeta_F$. We define the following functions on $F_\infty$ (depending only on $X$) for $s\in \C$, $s\not=0,1$ :

\vspace{2mm}

{\madef {For $x\in F_\infty$, $X=|x|$,
$$
H_s(x) = X^{s-1} \sum_{N\fa \le X} N\fa^{-s} - \frac{\kappa}{1-s}
$$
where $\fa$ runs over non zero \rm{integral} ideals  $\fa \subset \Oo$.}
$$
A_s (X) = \sum_{N\fa\le X} N\fa^{-s}.
$$}

\vspace{2mm}

{\madef For $x\in F_\infty$, $X=|x|$,

$$
K_s(x) = D^{-1/2}X^{s-1} \sum_{\scriptstyle \fa \subset \Dd^{-1} \atop \scriptstyle N\fa \leq X}  N\fa^{-s}-\kappa \frac{D^{1/2}}{1-s}.
$$}

\vspace{2mm}

These are functions of slow growth on $F_\infty$, and therefore tempered distributions. We call $H_s$ the \textit{Tate kernel}.

\vspace{2mm}

\setcounter{equation}{0}

{\montheo Assume $\sigma=Re(s) \in ]0,1[$. Then $\zeta_F(s)=0$ if, and only if 
$$
\Ff(H_s) = -K_{1-s}.
$$}

\subsection{}
The proof follows easily from Tate's functional equation. Let $\A$ denote the ad\`eles of $F$, $\A=F_\infty \times \prod\limits_v{}' F_v$ (restricted product), $v$ running over the finite places. We write $\A=F_\infty\times\A_f$. Let $I$ be the group of id\`eles, $I=F_\infty^\times \times I_f$, and $I^+$ the set of elements in $I$ of id\`{e}le norm $>1$. We endow $\A$ with the self--dual measure of Tate \cite{Ta}, associated to Tate's construction of local additive characters. Let $h=h_\infty \otimes h_f$, where $h_\infty \in \CSs(F_\infty)$, and $h_f=\bigotimes h_v$, $h_v=ch(\Oo_v)$, $\Oo_v\subset F_v$ being the integers. Let $\Dd_v$ be the local different in $F_v$, $\Dd_v^{-1}$ its inverse. Then\begin{equation}
\hat{h}_v= |\Dd_v|^{-1/2} ch(\Dd_v^{-1})
\end{equation}
where $ch(\ )$ denotes the characteristic function. The self--dual measures $dx_v$ define as usual Haar measures $d^{\times} x_v$, on $F_v^\times$, hence a measure on $I$ (at the Archimedean places, $d^\times x_v = \frac{dx_v}{|x_v|}$). With $\hat{h}= \hat{h}_\infty \otimes \bigotimes\limits_v \hat{h}_v$, Tate's formula reads 
$$
\begin{array}{rl}
Z(h,s) &=\displaystyle \int_{I^+}h(x) |x|^s d^\times x + \int_{I^+} \hat{h}(x)|x|^{1-s}d^\times x\\	
\noalign{\vskip2mm}
&-\dfrac{\kappa}{1-s} \hat{h}(0) -\dfrac{\kappa}{s}h(0).		
\end{array}
$$
Here $Z(h,s) = Z(h_\infty,s) D^{-1/2} \zeta_F(s)$, $Z(h_\infty,s)$ being the product of the local integrals $Z(h_v,s)$ $(v\mid \infty)$. Assume then that $\sigma \in ]0,1[$. The local, Archimedean, integrals, are holomorphic for $\sigma>0$ ; for a suitable choice of the $h_v$, they are non--zero. Therefore $s$ is a zero if, and only if, for all values of $h_\infty$ :
\begin{equation}
\int_{I^+} h(x) |x|^s d^\times x - \frac{\kappa}{1-s} \hat{h}(0) + \int_{I^+} \hat{h}(x)|x|^{1-s} d^\times x -\frac{\kappa}{s}h(0 )=0.
\end{equation}
In this formula, we may replace the finite part of the measure $d^\times x$ by the ``trivial'' measure on $I_f$, giving mass $1$ to the idelic units, $\prod\limits_v \Oo_v^\times := \Oo_f^\times$. We denote by $\gd^\times x$ this new measure on $I_f$ (and $I$, the Archimedean components being unchanged.) Then $d^\times x= D^{-1/2}\gd^\times x$, cf. Tate \cite[p.310]{Ta}. Therefore (1.2) yields :
\begin{equation}
\int_{I^+}h(x) |x|^s \gd^\times x - \dfrac{\kappa D^{1/2}}{1-s} \hat{h}(0) + \int_{I^+}\hat{h}(x)|x|^{1-s}\gd^\times x - \dfrac{ \kappa D^{1/2}}{s} h(0) =0.
\end{equation}
We now write the first part in classical terms. By (1.1) $D^{1/2}\hat{h}_f(0)=1$. We have
$$
\int_{I^+} h(x)|x|^s \gd^\times x = \int_{F_\infty} h_\infty(x_\infty)|x_\infty|^{s-1} dx_\infty \int_{I_f}^X h_f(x_f)|x_f|^s\gd^\times x_f
$$
where the second integral runs over
$$
\{x_f\in I_f \mid |x_f|\ge X^{-1}\}.
$$
The integrand is invariant by $\Oo_f^\times$. We have
\begin{equation}
I_f = \coprod_{\alpha=(\alpha_ v)} \Big(\prod_v \varpi_v^{\alpha_v}\Oo_v^\times)
\end{equation}
where the $\alpha_v$ are almost all $0$. They must be positive since $h_v$ is supported on the integers.
If $x_f$ is in the $\alpha$--component of (1.4), $|x_f| = \prod\limits_v q_v^{-\alpha_v}$. Thus $\prod q_v^{\alpha_v} \le X$, and the corresponding value of $|x_f|^s$ is $\prod q_v^{-\alpha_vs}$. If $\fa = \prod 
\pp_v^{\alpha_v}$, we see that the inner integral is
$$
\sum_{N\fa\le X} N\fa^{-s}.
$$

Taking into account the residue term in (1.3), we see that the first two terms in (1.3) yield
\begin{equation}
\int_{F_\infty} h_\infty (x_\infty) H_s(x_\infty) dx_\infty.
\end{equation}

We compute similarly the second part of (1.3). The integral is
$$
\int_{F_\infty} \hh_\infty(x_\infty)|x_\infty|^{-s} dx_\infty \int_{I_f}^X \hh_f(x_f) |x_f|^{1-s} \gd^\times x_f
$$
with identical notation. Now in the decomposition (1.4) only terms such that $\varpi_v^{\alpha_v} \in \Dd_v^{-1}$ occur : thus $\alpha_v \ge -\delta_v$, where $\delta_v =val(\Dd_v)$. Again $\prod q_v^{\alpha_v} \le X$, and the value of $|x_f|^{1-s}$ is $\prod q_v^{\alpha_v(s-1)}$. The inner integral is therefore --- since $\hh_f = D^{-1/2} ch(\Dd_f^{-1})$, $\Dd_f^{-1}= \prod\limits_v \Dd_v^{-1}$ by (1.1) --- 
$$
D^{-1/2} \sum_{\fa \subset \Dd^{-1}\atop N\fa \le X} N\fa^{s-1}.
$$
Finally, the second part of (1.3) is equal to the integral of $\hh_\infty$ against
$$
D^{-1/2} X^{-s} \sum_{\fa \subset \Dd^{-1}\atop N\fa \le X} N\fa^{s-1} - \frac{D^{1/2}\kappa}{s}.
$$

The equality (1.3) then implies Theorem 1.1. 

\subsection{}

We will make no analytic use of Theorem 1.1. We note, however, that a direct ``classical'' proof that $\zeta(s)=0$ implies the identity of Fourier transforms  is easily obtained for $F=\Q$, where $H_s=K_s$. In this case $H_s$ is the even function defined for $x\ge 0$ by
\begin{equation}
H_s(x) = x^{s-1} \sum_{n\le x} n^{-s} - \frac{1}{1-s}
\end{equation}
$(n\in \N-\{0\})$. Proceeding as in Titchmarsh \cite[p.14]{T1} we get another expression of $H_s$. By the Euler--Maclaurin formula, we have for $0<x<y$ and $s\in \C$ :
$$
\sum_{x<n\le y} n^{-s} = \int_x^y t^{-s} dt - \psi(y) y^{-s} + \psi(x) x^{-s} -s \int_x^y \psi(t) t^{-s-1}dt,
$$
with $\psi(x) = \{x\}-1/2$ $(x\ge 0)$.

For $\sigma>1$ we get
$$
\sum_{x<n} n^{-s} =- \frac{x^{1-s}}{1-s} + \psi(x)x^{-1} -s \int_x^\infty \psi(t) t^{-s-s-1}dt
$$
so
$$
\sum_{n\le x} n^{-s} = \zeta(s) +    \frac{x^{1-s}}{1-s} - \psi(x)x^{-s} +s \int_x^\infty \psi(t) t^{-s-1}dt.
$$
This is true for $\sigma>0$ $(s\not=1)$ by analytic continuation. If $\zeta(s)=0$, $0<\sigma<1$ we get 
\begin{equation}
H_s(x) = -\psi(x) x^{-1} +s x^{s-1} \int_x^\infty \psi(t)t^{-s-1}dt.
\end{equation}
Since the average value of $\psi$ vanishes, the last term of (1.7) is an $O(x^{-2})$. In  particular $H_s \in L^2(\R)$.

On the other hand, (1.6) yields the differential equation
\begin{equation}
\begin{array}{c}
D H_s =(s-1/2) H_s +R,\\
R=-1 +\displaystyle\sum_n{}' \delta_n
\end{array}
\end{equation}

where $D=x \frac{d}{dx} +\frac{1}{2}$, $\delta_n$ is the Dirac measure at $n\in \Z$, $\Sigma'$ runs over $(n\not=0)$, and $H_s$ is seen as a tempered distribution.

Taking the Fourier transform yields

$$
D\hat{H_s} = (1/2-s) \hat{H_s} -R 
$$
whence
$$
D(\hat{H_s}+H_{1-s}) = (1/2-s) (\hat{H_s}+H_{1-s}) 
$$

This implies that $F= \hat{H_s}+H_{1-s} = C x^{-s}$, say on $\R_+$, which is impossible, since $F$ is $L^2$, unless $C=0$. Thus

$$
\Ff H_s= -H_{1-s}.
$$

\subsection{}
We now return to the general case and consider the Mellin transform of $H_s$, viewed as a function of $X\in \R_+$. In order to define
\begin{equation}
\Mm H_s(w) = \int_0^\infty H_s(x) x^{w-1}dx
\end{equation}
(where we write $x\ge 0$ for $X$) we must control the order of growth of $H_s$ at infinity.

However, we first consider (1.9) for $F=\Q$ where the calculation is explicit. From (1.7) we see that the integral is absolutely convergent for $0<\tau<1$, where $\tau=Re(w)$. Using (1.6) we consider
\begin{equation}
\int_0^X \Big(x^{s-1} \sum_{n\le x} n^{-s} -\frac{1}{1-s}\Big) x^{w-1}dx.
\end{equation}
There is a first term equal to
$$
(a)\ \ - \frac{1}{1-s} \int_0^X x^{w-1} dx = -\frac{1}{w(1-s)}X^w.
$$
The other term is
$$
\sum_{n\le X} n^{-s} \Big[\frac{x^{s+w-1}}{s+w-1} \Big]_n^X = \frac{1}{s+w-1} \sum_{n\le X} \{n^{-s} X^{s+w-1} - n^{w-1}\}.
$$
By the formulas in \S 1.3,
$$
\sum_{n\le X} n^{w-1} = \zeta(1-w) + \frac{X^w}{w} +O(X^{\tau-1} ).
$$

On the other hand, since $\zeta(s)=0$,
$$
\sum_{n\le X} n^{-s} = \frac{X^{1-s}}{1-s} +O(X^{-\sigma}),
$$
so
$$
X^{s+w-1} \sum_{n\le X} n^{-s}= \frac{1}{1-s} X^w +O(X^{\tau-1}
).
$$
The other term of the integral (1.10) is therefore equal to
$$
(b)\ \ -\frac{\zeta(1-w)}{s+w-1} + \frac{1}{s+w-1} \Big\{ -\frac{1}{w} + \frac{1}{1-s}\Big\} X^w +o(1).
$$
This yields

\vspace{2mm}

{\maprop
$(F=\Q)$. --- For $0<Re(w)<1$,
$$
\Mm H_s(w) = \frac{\zeta(1-w)}{1-s-w}.
$$}

\vspace{2mm}

We will now see that the same formula is true for any number field (but valid analytically in a smaller domain.) However we must proceed differently. To control the growth of $H_s$ for $x\rg \infty$ is to control the growth of $A_s$.

In this section we will simply derive an estimate from the estimate of $A_s(x)$ coming from a classical theorem of Landau \cite{La}\footnote{Cf Michel \cite[Thm 1.2]{Mi}. We return to Landau's theorem in \S~2.3 (cf. Theorem~2.1)}. We have for $x\ge 1$ :
$$
A_s(x) = \sum_{n\le x}a_nn^{-s}
$$
where $a_n = \sum\limits_{N\fa=n}1$. The z\^eta function
$$
L(s) =\zeta_F(s) = \sum a_n n^{-s} \quad (Re\ s > 1)
$$
can be seen as an $L$--function over $\Q$. For further reference we note that it is then an Euler product with factors
\begin{equation}
L_p(s) =\prod_{i\le m_p} (1-\zeta_ip^{-s})^{-1}
\end{equation}
with $m_p \le m=[F:\Q]$, and $m_p=m$ at the unramified primes, the $\zeta_i$ being roots of unity. Landau's assumptions are verified, and we get for
$$
A(x) =A_0(x) = \sum_{n\le x} a_n
$$
the estimate
\begin{equation}
A(x) = \kappa \ x+O(x^{\frac{m-1}{m+1}+\varepsilon}).
\end{equation}
By Stieltjes integration, we obtain the following estimate for $A_s(x) = \int_{1^{-}}^x\break t^{-s} dA$.

We fix $\varepsilon$ (small) and set $\mu=\frac{m-1}{m+1}+\varepsilon$. We assume $\sigma=Re (s) <1$. Then
$$
\begin{array}{rll}
A_s(x) &= \dfrac{\kappa}{1-s} x^{1-s} +O(x^{\mu-\sigma}) & (\sigma<\mu)\\
\noalign{\vskip2mm}
	&= \dfrac{\kappa}{1-s} x^{1-s} +O(1) &(\sigma>\mu).
\end{array}
$$
Thus
$$
\begin{array}{rll}
H_s(x)= &O(x^{\mu-1}) &(\sigma<\mu)\\
&O(x^{\sigma-1}) &(\sigma>\mu),
\end{array}
$$
and the integral (1.9) converges in the following domain, setting $\tau=Re (w)$~:
$$
\begin{array}{rl}
\tau <1-\mu &(\sigma<\mu)\\
\tau < 1-\sigma &(\sigma>\mu).
 \end{array}
$$
In particular we record that

\vspace{2mm}

\setcounter{monlem}{0}
{\monlem 

$(i)$ For $\sigma<\frac{m-1}{m+1}$, $\Mm H_s(w)$ is given by an absolutely convergent integral for $0<\tau <\frac{2}{m+1}$

$(ii)$ For $\sigma>\frac{m-1}{m+1}$, this is true  for
$$
0<\tau<1-\sigma.
$$}


So far we have not assumed that $s$ was a zero of $\zeta_F$. We now do so ; also $0<\sigma<1$. Perron's formula yields for $x>0$
$$
A_s^*(x) = \sum_{n\le x} \\ ^* a_n n^{-s} = \frac{1}{2i\pi} \int_{c-i\infty}^{c+i\infty} \zeta_F (s+w) \frac{x^w}{w} dw
$$
for $c>1-\sigma$. Here $\sum^*$ means as usual that the last term is pondered by $\frac{1}{2}$ if $x$ is an integer. The vertical integrals are Cauchy principal values. Thus
$$
x^{s-1}A_s^* (x) = \frac{1}{2i\pi} \int_c \zeta_F (s+w) \frac{x^{s+w-1}}{w} dw
$$
(where $\int_c$ will denote $\int_{c-i\infty}^{c+i\infty})$
$$
=\frac{1}{2i\pi} \int_{c'} \zeta_F (1-w) \frac{x^{-w}}{1-s-w} dw
$$
by the change of variable $w \mapsto 1-s-w$ ; here
$$
c'=1-\sigma-c<0.
$$

Consider now $c''>0$. We shift the integration to the abscissa $c''$ ; we assume $c''<1$. Then, formally,
$$
\int_{c'} = \int_{c''} - 2i\pi \sum Res
$$
where $\sum Res$ is the sum of the residues in the band bounded by the two lines. There is a residue equal to $\frac{-1}{1-s}\kappa$, at $w=0$ ; the denominator $\frac{1}{1-s-w}$ contributes a pole at $w=1-s$, but $\zeta_F(w)=0$. Therefore our last integral~is
$$
\frac{1}{2i\pi} \int_{c''}\zeta_F(1-w) \frac{x^{-w}}{1-s-w} + \frac{\kappa}{1-s}.
$$

To justify the translation of the integral, we must now consider (for $T\rg \infty$) the horizontal integrals on $[c'+iT,\ c''+iT]$. We have
$$
1-c'' \le Re(1-w) \le 1-c' \quad (c'<0).
$$
We have assumed $c''\in ]0,1[$. Then the convexity estimate, for $w=\tau+it$, $\tau\in [c',c'']$, is 
$$
|\zeta_F(1-w)| \ll t^{\frac{m}{2}c''+\varepsilon},
$$
cf. Iwaniec-Kowalski \cite[p.100, p.126]{IK}\footnote{The statement of Theorem 5.30 in \cite{IK} is incorrect.}. However it is known \cite[Thm.1.1]{MV} that \textit{subconvexity} holds for the z\^eta functions of number fields, so in fact
$$
|\zeta_F(1-w)| \ll t^{\frac{m}{2}c''-\delta}
$$
for some small $\delta>0$. Consequently the horizontal integrals tend to $0$ if
$$
c'' \le \frac{2}{m}.
$$

Therefore :

\vspace{2mm}

{\monlem Assume $c'' \le \frac{2}{m}$. Then, for $x>0$,
$$
H_s^*(x) = \frac{1}{2i\pi} \int_{c''} \frac{\zeta_F(1-w)x^{-w}}{1-s-w}dw.
$$}

In this paper we will be mostly concerned with the integral on the critical line $Re(w)=\frac{1}{2}$. We must therefore assume $m\le 4$. Furthermore, we will consider the quadratic integral $\int|\frac{\zeta_F(1/2-it)}{1/2 -it-s}|^2dt$ ; even with the (known) subconvex estimates, this converges only for $m\le 2$. We now assume this\footnote{We will return to these statements, using known conjectures, in \S 3.4.}. As a direct consequence of Lemma 1.2, we now have~:

\vspace{2mm}

{\maprop Assume that $F$ is $\Q$ or a quadratic extension of $\Q$. Let $s$ a zero of $\zeta_F$ with $0<\sigma<1$. Then

$(i)$ $H_s \in L^2 (\R_+,dx)$

$(ii)$ $\Mm H_s(w) = \frac{\zeta_F(1-w)}{1-s-w}$ $(Re\ w=1/2)$.}

\vspace{2mm}

In part $(ii)$, $\Mm H_s(w)$ is a priori defined as an $L^2$ function of $w$ ; for $w=\frac{1}{2} +it$,
$$
\Mm H_s(w) = \int_0^\infty H_s(x) x^{1/2+it} \frac{dx}{x}.
$$
i.e.
\begin{equation}
\Mm H_s(w) = \int_{-\infty}^\infty (H_s(e^X) e^{X/2}) e^{itX}dX,
\end{equation}
and $(i)$ is equivalent to the fact that $H_s(e^X) e^{X/2} \in L^2(\R,dX)$. Cf. Titchmarsh \cite[Thm. 7.1]{T2}.

The proof is clear : write
\begin{equation}
H_s^*(x) = \frac{1}{2i\pi} \int_{-\infty}^{+\infty} \frac{\zeta_F(1/2-it)x^{-1/2-it}}{\frac{1}{2}-it-s}\  d(it)
\end{equation}
according to Lemma 1.2. On the other hand, in view of the last remarks, we can see the right--hand side as an additive Fourier transform (evaluated at $X=\log x)$. Thus the Fourier transform (in $L^2$) of the right--hand side of (1.14) coincides with $e^{X/2}H_s^*(e^X)$. This implies $(i)$, and $(ii)$ by the involutivity  of the Fourier transform.

It is interesting to notice the relation with the constraints in Lemma 1.1. Assume $m=2$. If $\sigma<1/3$, $\Mm H_s(w)$ is holomorphic for $\tau<\frac{2}{3}$. If $\sigma>\frac{1}{3}$, it is holomorphic for $\tau<1-\sigma$. In particular, $\tau=1/2$ is attained only if $\sigma<\frac{1}{2}$. For $\sigma=1/2$, we are on the boundary of the domain of convergence\footnote{But recall that we are using only the crude estimates given by Stieljes integration. Compare Theorem 2.2 and its Corollary.}.

\vspace{2mm}

\noindent\textbf{{Corollary.}} \textit{With the assumptions of Proposition $1.2$,}
$$
\int_0^\infty |H_s(x)|^2 dx = \frac{1}{2\pi} \int_{-\infty}^{+\infty} \Big|\frac{\zeta_F(\frac{1}{2}+it)}{1/2+it-s}\Big|^2 dt.
$$

This follows from $\int_0^\infty|H_s(x)|^2 dx=\int_{-\infty}^\infty |H_s(e^X)e^{X/2}|^2dX$, cf \cite{T2}.

We can now prove the main theorem of this chapter :

\vspace{2mm}

{\montheo Assume $s$ is a zero of $\zeta_F$ in the critical strip, and consider
$$
I(s) =\int_{-\infty}^{+\infty} \Big|\frac{\zeta_F(\frac{1}{2}+it)}{1/2+it-s}\Big|^2 dt.
$$
Then
$$
I(s) \ge 2\pi \Big(\log 2 - \frac{2\kappa}{|1-s|}\Big)
$$
}

We recall that $\kappa$ can be majorized, for any $\varepsilon>0$, by $C(\varepsilon)D^\varepsilon$ where the constant is explicit (cf. e.g. Lang \cite[Ch.XVI]{Lang}). In particular the lower bound is effective.

The proof is now obvious : the integral of $|H_s|^2$ is larger than
$$
\int_1^2 |H_s(x)|^2 = \int _1^2 \Big|x^{s-1}-\frac{\kappa}{1-s}\Big|^2 dx,
$$
which is bounded below by
$$
\int_1^2 \Big(x^{2\sigma-2}-2 \kappa \frac{ x^{\sigma-1}}{|1-s|}\Big) dx.
$$

The first term yields
\begin{eqnarray}
&\log 2 &\mathrm{for} \ \sigma=1/2\\
&f(\sigma) =\dfrac{2^{2\sigma-1}-1}{2\sigma-1} &\mathrm{for}\ \sigma\not=1/2.
\end{eqnarray}

For $\sigma>\frac{1}{2}$, the value of (1.16) is larger than $\log 2$. However $I(s)=I(1-s)$, as follows from $\zeta_F(s)= \zeta_F( \bar{s}).$

(If $\sigma= Re(s)$ is sufficiently close to 1, one can obtain better bounds by considering the integrals on the intervals $[N,N+1]$, $N\le 3$, using that $\frac{1}{2^\sigma}+\frac{1}{3^\sigma}<1$.)

\section{Estimates for $H_s$ : general case}

\subsection{} In this chapter we return to the problem, raised in \S~1.4, of getting estimates for $H_s$ better than that given by Landau's theorem and partial integration. However we find that the question can be posed --- and partly solved, giving non trivial estimates --- in  a very general framework.

We will consider the function $H_s$ associated to very general Dirichlet series. For definiteness assume first that $\pi$ is a cuspidal, unitary representation of $GL(m,\A_\Q)$. Thus
$$
\pi = \pi_\infty \otimes \bigotimes_p \pi_p.
$$
Its $L$--function $L(s,\pi)$ is an Euler product
$$
L(s,\pi) = \prod_p L(s,\pi_p) \qquad (Res>1)
$$
convergent for $Res>1$ \cite{JS}.
We assume first that $\pi$ verifies the Ramanujan conjecture at all primes, i.e., that $\pi_p$ and $\pi_\infty$ are tempered. After the early work of Deligne and Serre, this is now known in many non--trivial cases, cf \cite{Ca}, \cite{Cl3}. The $L$--function can then be written
$$
L(s,\pi) = \sum_1^\infty a_n \ n^{-s}
$$
with
\begin{equation}
|a_n| \le \tau_m(n) \overset{}{\underset{\varepsilon}{\ll}}  n^\varepsilon
\end{equation}
as in \S 1.4. (In particular the estimate is uniform for $m$ and $\varepsilon$ fixed.) There is a well--defined $L$--function $L(s,\pi_\infty)$, which can be written \footnote{In the next formula, and others, the reader will not confuse $\pi = 3.14159...$ and $\pi$ = automorphic representation. }
\begin{equation}
L(s,\pi_\infty) =c(\pi_{\infty}) \pi^{-\frac{ms}{2}} \prod_{i=1}^m \Gamma\big(\frac{s+c_i}{2}\big).
\end{equation}
This will be reviewed presently. There is also a conductor $D$, a positive integer. See Jacquet \cite{J1}, as well as \cite{J2}.

We write
$$
\Lambda(s,\pi)= D^{s/2} L(s,\pi_\infty)L(s,\pi).
$$
Then $\Lambda(s,\pi)$ satisfies a functional equation :
\begin{equation}
\Lambda
(s,\pi) = \varepsilon(\pi) \Lambda(1-s,\tilde{\pi}).\end{equation}
where $|\varepsilon(\pi)\vert=1$ and $\tilde{\pi}$ is the dual representation. It is holomorphic in the whole plane, unless $m=1$ and $\pi(x)=|x|^a$ $(x\in \A_\Q^\times)$ with $a\in i\ \R$.

More generally, we can consider $\pi = \pi_1 \times \cdots \times \pi_r$, $\pi_i$ being cuspidal tempered representations of $GL(m_i,\A_\Q)$. We have similar properties for $L(s,\pi)$ ; its poles are easily described.

There are other $L$--functions with the same properties. First let $F/\Q$ be a finite extension (not necessarily Galois) and let $\pi_F$ be a cuspidal representation of $GL(\mu,\A_F)$. By the Langlands classification (now known at all primes, cf.\cite{HT}) $\pi_F$ defines a representation
$$
\pi_\Q = \Ind_F^\Q \pi_F
$$
of $GL(n,\A_\Q)$, of degree $m=\mu d$, $d=[F:\Q]$. Unless $F$ is a soluble Galois extension of $\Q$, $\pi_\Q$ is not known to be automorphic. We assume that $\pi_F$ verifies the Ramanujan conjecture (at all primes). Then so does $\pi_\Q$ ; when base change is known  \cite{AC} it is of the form $\pi_1 \times\cdots\times \pi_r$ as before. (Note that the local correspondence has to be normalised, as can be done, so induction on the (local) Weil group side preserves tempered representations). We have, for $Re(s)>1$,
$$
L(s,\pi_F) = L(s,\pi_\Q)
$$
by the inductivity of $L$--functions ; both sides can be seen as Euler products over $\Q$. Defining the Archimedean factor $L(s,\pi_{\Q,\infty})$ as before, we see that $L(s,\pi_\Q)$ extends meromorphically to the complex plane ; we can form $\Lambda(s,\pi_\Q)$ with the same properties as before.

If $\mu=1$ and $\pi_F$ is trivial $L(s,\pi_F)=\zeta_F(s)$ so this accounts for the situation in Chapter 1. (For $m\geq2$, in the cases where the Ramanujan conjecture is known, $F$ will be totally real, or a $CM$ field.)

Similarly, assume $\rho_F$ is an irreducible Artin representation of $\Gal(\overline{F}/F)$ such that the Artin Conjecture is known for the Artin $L$--function $L(s,\rho_F)$\footnote{I do not know if the Artin Conjecture is known for some $\rho_F$ such that the associated representation $\pi_F$ has not been shown to exist.}. Again we can consider $\rho_\Q = \ind _{G_F}^{G_\Q}\rho_F$. Its $L$--function is just $L(s,\rho_F)$ (seen as an Euler product over $\Q$.) It has the specified properties, even if the representation $\pi_\Q$ associated to $\rho_\Q$ is not known to exist.

\subsection{} In all these cases we consider
\begin{equation}
L(s) = L(s,\pi_\Q) = \sum_n a_n n^{-s},
\end{equation}
defined for $Re(s)>1$, and the Archimedean factor $L(s,\pi_\infty)$. The  representation $\pi_\infty$ is well--defined and tempered.

We now describe the factor $L(s,\pi_\infty)$. By the Langlands classification\footnote{Again, this is Langlands's normalisation, compatible with unitarity and induction, cf. \cite{Cl1}.}, $\pi_\infty$ is associated to a representation $r(\pi): W_\R \rg GL(m,\C)$ where $W_\R$ is the real Weil group. We have
$m=m_1+2m_2$, $r=r(1)+r(2)$, $r(1)$ is a sum of $m_1$ $1$--dimensional representations associated to characters $\nu_i$ of $\R^\times$, and $r(2)$ is a sum of $m_2$ $2$--dimensional representations
$$
r_j =\ind_{W_\C}^{W_\R}\ \mu_j:= r(\mu_j)
$$
where $\mu_j$ is a unitary character of $\C^\times$ such that $\mu_j (\bar{z}) \not  \equiv \mu_j(z)$.

A unitary character $\nu$ of $\R^\times$ can be written
$$
\nu(x)=(\mathrm{sgn} \ x)^\varepsilon\ |x|^c\quad (\varepsilon=0,1;\ c\in i \R).
$$
Its $L$--function is
\begin{eqnarray}
L(s,\nu) &=& \pi^{-\frac{s+c}{2}}\ \Gamma\big(\dfrac{s+c}{2}\big)\quad (\varepsilon=0)\\
L(s,\nu) &=&\pi^{-1/2-\frac{s+c}{2}}\ \Gamma\big(\dfrac{s+c+1}{2}\big)\quad (\varepsilon=1).
\end{eqnarray}
A unitary character $\mu$ of $\C^\times$ such that  $\mu_j (\bar{z}) \not  \equiv \mu_j(z)$ can be written
$$
\mu(z) =z^p\overline{z}^q,\ p\not=q,\ p-q\in \Z,\ p+q\in i\R.
$$
Its $L$--function, equal to that of $r(\mu)$, is
\begin{eqnarray}
L(s,\mu) &=& 2(2\pi)^{-q-s}\ \Gamma(q+s) \quad (p-q<0)\\
L(s,\mu) &=& 2(2\pi)^{-p-s}\ \Gamma(p+s) \quad (p-q>0).
\end{eqnarray}
(see Tate \cite[\S~3]{Ta2}.)

Using the duplication formula, we see that the factor (2.7) is equal to
\begin{equation}
L(s,\mu) = \pi^{-1/2}\ \pi^{-s-q}\ \Gamma\big(\frac{s+q}{2}\big) \Gamma\big(\frac{s+q+1}{2}\big).
\end{equation}
It is then easy to see that $L(s,\pi_\infty)$ is given by an equality (2.2). We now rewrite the functional equation (2.3) as
$$
L(1-s,\pi) = \varepsilon(\pi)\gamma(s)\ L(s,\tilde{\pi})
$$
with
\begin{equation}
\gamma(s) = (\pi^{-m}D)^{s-1/2}\ \frac{c(\tilde{\pi}_\infty)\Gamma(s,\tilde{\pi}_\infty)}{c(\pi_\infty)\Gamma(1-s,\pi_\infty)},
\end{equation}
$\Gamma(s,\pi_\infty)$ being the product of the gamma factors.

\subsection{} Given an $L$--function $L(s) =L(s,\pi)$, as in \S~2.1, we now consider for $x\ge0$
\begin{eqnarray}
A_s(x) &=& \displaystyle\sum_{n\le x} a_n \ n^{-s}\\
H_s(x) &=& x^{s-1} \displaystyle\sum_{n\le x} a_n\ n^{-s} -\dfrac{\kappa}{1-s}
\end{eqnarray}
with $\kappa=Res_{s=1}L(s)$. We will assume that $L(s)$ has at most a simple pole at $s=1$. In many cases, Landau (in 1915 !) had already obtained an excellent estimate for $A_0(x) = \sum\limits_{n\le x}a_n.$ It is straightforward to check that his conditions (I--VII), p.~2--4 of \cite{La}, are satisfied. However, he assumes that the $\Gamma$--factors, $\Gamma(\frac{s+c_i}{2})$, are given by real parameters $c_i$. This means here that the real characters $\nu_i$ are equal to 1 or to $
\rm sgn( x)$ ; in the complex case, that $\mu(z)=z^p\ \overline{z}^q$ with $p,q$ real, i.e. that
$$
\mu(z) = (z/\overline{z})^p, \qquad p\in \frac{1}{2} \ \Z.
$$
Landau does not assume that the $a_n$ are real, but that $|a_n|\le d_n$ where
$$
L^*(s) = \sum_{n=1}^\infty d_n\ n^{-s}
$$
has a functional equation of the same type as that of $L(s)$. Since our representation $\pi$ (or the $\pi_i$) is tempered at all primes, $L(s)$ is given by an Euler product
$$
L(s) = \prod_p L_p(s) \qquad (Res>1)
$$
with factors of type (1.11), the $\zeta_i$ being now complex numbers with absolute value $1$, at unramified primes. Assume $p$ is a ramified prime. One can write $\pi_p= \pi_1 \times \pi_2 \times... \times \pi_r$, the local  analogue of the global decomposition of $\pi$ in $\S$ 2.1. Each $\pi_j$ is a unitary discrete series representation of $GL(m_j , \Q_p)$. Then $L(\pi,s)= \prod L(\pi_j,s)$, the product running over the $j$ such that $\pi_j$ is a unitary twist of the Steinberg representation. The corresponding factor is 
$$
(1-p^{\frac{1-m_j}{2}+i t_j} p^{-s})^ {-1}
$$
with $t_j$ real. Cf. \cite{J1}, \cite {GJ} \footnote{However Jacquet's expression in \cite{J1} seems incorrect. Compare \cite{GJ}, Theorem 7.11.}. It follows that we can take $L^*(s)=\zeta_\Q(s)^m$. The following result was therefore known in 1915 (for $m>1$, which is equivalent to condition (9) in \cite{La})~:

\vspace{2mm}

{\montheo (Landau) Assume that $m>1$ and that each real character of $\R^\times$ occuring in $r(\pi)$ is of finite order, and each complex character of the form $(z/\overline{z})^p$, $p\in \frac{1}{2}\Z$. Then
$$
A_0(x) = \kappa\ x+O(x^{\frac{m-1}{m+1}+\varepsilon}).
$$}

In particular, we see that Lemma 1.1 remains true under these assumptions. It is remarkable that this applies to the cuspidal \textit{algebraic} representations of \cite{Cl1} (twisted so as to render them unitary) when they satisfy the Ramanujan hypothesis. However, we want to do better.

\subsection{} In order to obtain an estimate of $A_s(x)$ better than that obtained by Stieljes integration, we now use the method of Friedlander--Iwaniec \cite{FI}.

Return to the equation (2.10). Friedlander--Iwaniec assume that $\gamma(s)$ is a product 
\begin{equation}
(\pi^{-m} D)^{s-1/2} \prod_{i=1}^m\ \frac{\Gamma(\frac{s+c_i}{2})}{\Gamma(\frac{1-s+c_i}{2})}.
\end{equation}
For the real factors (2.5), (2.6), the parameter $c$ is sent to $-c$ when $\pi$ is sent to $\tilde{\pi}$ ; similarly for (2.7), (2.8), $(p,q)$ is sent to $(-p,-q)$.  We therefore now make the following assumption\footnote{This assumption is also (implicitly) made by Iwaniec and Kowalski in \cite{IK} : see p. 94, after (5.3), and the equality $\gamma(\overline{f},s)=\gamma(f,s)$, p.~94,~l.~-9.}

\vspace{2mm}

{\myhypo $\pi_\infty$ is self--dual.}

\vspace{2mm}

In this case $\gamma(s)$ is given by a product (2.13) with $c_i\in \C$ and furthermore
\begin{equation}
Re(c_i)\ge 0.
\end{equation}
(For the parameters $c_i$ coming from (2.7), the condition on $Re(c_i)$ follows from $Re(p+q)=0$ and $p-q<0$ ; \textit{idem} for (2.8.)

\vspace{2mm}

\noindent\textbf{Remark.} If $\pi$ is associated to a Galois representation of $\Gal(\overline{F}/F)$, $\pi_\infty$ is self--dual.

The expression (2.13) shows that $\gamma(s)$ is holomorphic for $\sigma>0$. Write $Q=D^{1/m}/2\pi$, so $2\pi Q \ge 1$.

\vspace{2mm}

{\monlem There exist constants $\omega_1,\omega_2 \in \C$ such that for $\sigma\ge \frac{1}{2}$, $t\ge 1$
\begin{eqnarray*}
\gamma(\sigma+it) &=& \omega_1 (Qt)^{m(\sigma-1/2)}(Qt/e)^{imt}(1+O(t^{-1}))\\
\gamma(\sigma-it) &=& \omega_2 (Qt)^{m(\sigma-1/2)}(Qt/e)^{-imt}(1+O(t^{-1})),
\end{eqnarray*}
the implicit constant in the $O$--terms being uniform (for $\pi$ given) for $1/2\le \sigma\le A$, $A>\frac{1}{2}$.}

\vspace{2mm}

(In fact the uniformity is obtained if the real parameters $(c_i)$ are fixed.)

This is essentially (1.6) in \cite{FI}, except that we do not obtain the same constants $\omega_i$. We sketch the proof. For $t\ge 1$, we have by Stirling's formula~\cite{Er}
\begin{eqnarray}
\Gamma
(\sigma+it) &=& \sqrt{2\pi}\ e^{(\sigma-1/2)i\frac{\pi}{2}}\ t^{\sigma-1/2}\ e^{-\frac{\pi}{2}t}(t/e)^{it}(1+O(1/t))\\
\Gamma
(\sigma-it) &=& \sqrt{2\pi}\ e^{(1/2-\sigma)i\frac{\pi}{2}}\ t^{\sigma-1/2}\ e^{-\frac{\pi}{2}t}(t/e)^{-it}(1+O(1/t))
\end{eqnarray}
It suffices to consider one of the quotients in the formula (2.13), associated to $c\in \C$, with $Re(c)\ge 0$. Let $q_1,\ q_2,\ q_3,\ q_4$ be the quotient associated to each non--trivial factor in (2.15), (2.16).

Let $c=\alpha+i\tau$. We find for example
\begin{eqnarray*}
q_1 = e^{i(\frac{\sigma+\alpha}{2}-\frac{1}{2})\frac{\sigma}{2}} &:&e^{-i(\frac{1-\sigma+\alpha}{2}-\frac{1}{2})\frac{\pi}{2}}\\
&=& e^{i(\alpha-1/2)\frac{\pi}{2}}.
\end{eqnarray*}
One computes similarly $q_2=(t/2)^{\sigma-1/2}(1+O(t^{-1}))$, $q_3= e^{-(\pi/2)\tau}$, $q_4=(t/2e)^{it}(1+O(t^{-1}))$. Multiplying all the factors, one finds for the product of quotients of $\Gamma$--functions :
$$
q(t) = e^{-im\frac{\pi}{4}}\ e^{i\frac{\pi}{2}c} (t/2e)^{imt}(t/2)^{m(\sigma-\frac{1}{2})} (1+O(t^{-1}))
$$
where $c$ is now equal to $\Sigma c_i$.

Multiplying by the factor $(\pi^{-m}D)^{s-1/2}$, one obtains the expression for $t>0$ $(\ge 1)$, with
$$
\omega_1=e^{-m\frac{i\pi}{4}} \ e^{i\frac{\pi}{2}c}.
$$
For $t<0$, using that, with obvious notation :
$$
\gamma(\overline{s},c_i) = \overline{\gamma(s,\bar{c}_i)}
$$
one obtains the requested expression, with $\omega_2=e^{im\frac{\pi}{4}}\ e^{-\frac{i\pi}{2}c}$. If the $c_i$ are real, we obtain the result of \cite{FI}. This will be the case for the algebraic representations of \cite{Cl1}, which includes the only cases where the Ramanujan conjecture is known.

\vspace{2mm}

\noindent\textbf{Remark.} The function $\gamma$ does not have the same expression if $\pi_\infty$ is not self--dual. Assume for instance that a factor (2.5) occurs with $c\not=0$, while the factor associated with $(-c)$ does not occur. Computing in the same fashion, we find for the factor $\frac{\Gamma(\frac{s+c}{2})}{\Gamma(\frac{1-s-c}{2})}$ an expression of the form $q_1\ q_2\ q_3\ q_4$, with $q_1 = e^{-\frac{i\pi}{4}}$, $q_2= \big(\frac{t}{2}\big)^{\alpha+\sigma-1/2}(1+O(t^{-1}))$, $q_3=1$ and $q_4=(t/2e)^{it}\ e^{i\tau}(t/2e)^{i\tau}(1+O(t^{-1}))$. The complete quotient is therefore of the form 
$$
q(t)=e^{-i\pi/4}\ (t/2e)^{i(t+\tau)}(t/2)^{\alpha+\sigma-1/2}\ e^{i\tau}\ (1+O(t^{-1})).
$$
We do not know if the argument of Friedlander--Iwaniec extends to this case.

\subsection{} We now imitate the proof of \cite{FI} in order to obtain an estimate for $A_s(x)$. \textit{We will assume} that $0<\sigma<1$, the case of $s=0$ being treated in \cite{FI}. We may further assume that $|x-N| \ge \frac{1}{4}$, $N$ being  the integer closest to~$x$.

Fix $\varepsilon>0$ (small) and let $c=1-\sigma +\varepsilon$. For $ 1 \leq T \leq x$, we have according to \cite[Lemma 3.12]{T1} :
$$
A_s(x) = \frac{1}{2i\pi} \int_{c-iT}^{c+iT} \ L(s+w) \frac{x^w}{w} dw + O\Big(\frac{x^c}{T\varepsilon^m} \Big) + O\Big(\frac{x^{1-\sigma+\varepsilon'}}{T}\Big)
$$
for any sufficiently small $\varepsilon'$. The implicit constants depend only on $\varepsilon$ and $\varepsilon'$. If $s=\sigma+i\tau$, we obtain by a change of variables
$$
A_s(x) = \frac{1}{2i\pi} \int_{\alpha+i(\tau-T)}^{\alpha+i(T+\tau)} \ L(w) \frac{x^{w-s}}{w-s} dw + O_\varepsilon\Big(\frac{x^{1-\sigma+\varepsilon}}{T} \Big).
$$
We may assume $\tau\ge 0$ : the computation for $\tau \leq 0$ is obviously similar. Here $\alpha=1+\varepsilon$. We shift the integral to the line $Re(w)=-\varepsilon$. We pass poles at $w=1$, with residue $\frac{\kappa}{1-s}x^{1-s}$, and at $w=s$, with residue $L(s)$. We write $w=u+it$. We now have to evaluate the integrals in the segments $Im(w)= \tau\pm T$, $Re\ w\in[-\varepsilon,1+\varepsilon]$.

We first do so, using only the convexity estimate. For $Re(w)=1+\varepsilon$,
$$
|L(w)| \le \zeta_\Q (1+\varepsilon)^m \ll \varepsilon^{-m} \ll 1.
$$
For $Re(w)=-\varepsilon$, $t\ge 1$
$$
L(w) \ll (Qt)^{m(1/2+\varepsilon)}
$$
by the functional equation. If $R=(Qt)^{m(1/2+\varepsilon)}\ge 1$, we have for\break $u \in [-\varepsilon,1+\varepsilon]$
$$
L(w) \ll R^{\frac{1+\varepsilon-u}{1+2\varepsilon}}
$$
whence
$$
\int_{-\varepsilon}^{1+\varepsilon} |L(u+it)| x^u du \ll x^{1+\varepsilon} + R\ x^{-\varepsilon} \ll x^{1+\varepsilon}+R.
$$
We will assume $T\ge 2\tau$, so the correction term is  $O\big(\frac{x^{-\sigma}}{T} \big(x^{1+\varepsilon} +\break(QT)^{m(1/2+\varepsilon)}\big)\big)$ and, $I_{-\varepsilon}$ being the integral for $Re(w)=-\varepsilon$ :
\begin{equation}
A_s(x) = I_{-\varepsilon} + \frac{\kappa}{1\!-\!s}\ x^{1-s}+ L(s) + O\Big(\frac{x^{1-\sigma+\varepsilon}}{T}\Big) + O \Big(\frac{x^{-\sigma}(QT)^{m(\frac{1}{2}+\varepsilon)}}{T}\Big).
\end{equation}
Applying the functional equation, we get
$$
I_{-\varepsilon} = \frac{\varepsilon(\pi)}{2i\pi} \int_{\alpha-i(T+\tau)}^{\alpha+i(T-\tau)} \gamma(w)L(w) \frac{x^{1-w-s}}{1-w-s} \ dw.
$$
The series for $\tilde{L}(w)=L(w,\tilde{\pi})$ being absolutely convergent, we get, with
$$
\tilde{L}(w) = \sum_{0}^\infty b_n\ n^{-w} :
$$
\begin{eqnarray}
I_{-\varepsilon} &=& \varepsilon(\pi) x^{1-s} \displaystyle\sum_0^\infty b_n\ c(nx),\\
c(y) &=& \dfrac{1}{2i\pi} \displaystyle \int_{\alpha-iT''}^{\alpha+iT'} \gamma(w) \dfrac{y^{-w}}{1-w-s}\ dw.
\end{eqnarray}
(Cf. (2.5) in \cite{FI}, for $s=0$.).

Recall that $T'=T-\tau$, $T''=T+\tau$ and that we assume
\begin{equation}
T\ge 2\tau.
\end{equation}
As in \cite{FI} we first consider $y >2(QT)^m$. Up to a constant, the  integral (2.19) is (Lemma 2.1)
$$
Q^{m(\alpha-1/2)} \int_{-T''}^{T'} \frac{|t|^{m(\alpha-1/2)}(Q|t|/e)^{imt}y^{-\alpha-it}}{-\varepsilon-it-s} \ (1+O(t^{-1})) dt.
$$

The function $F(t)=m\ t\ \log(Q|t|/e)-t\ \log y$ has derivative
$$
\begin{array}{ll}
F'(t) = m \log |Qt| - \log\ y. \\
\end{array}
$$
Since $y>2(QT)^m$, $F'(t)<-\log 2$ for $t\in [-T,T']$. Moreover, $F'$ is monotone on $(t>0)$ and $(t<0)$. By \cite[Thm~2, p.~104]{Te}, we deduce that $\int_{T_1}^{T_2}e^{iF(t)}dt \ll 1$ for $-T\le T_1\le T_2\le T'$. Let $G(t) =\int_0^t e^{iF(t)}dt$, which is thus absolutely bounded in the interval. The integral (2.19) has a first term equal to
$$
Q^{m(\alpha-1/2)} y^{-\alpha} \int_{-T}^{T'} \frac{|t|^{m(\alpha-1/2)}e^{iF(t)}}{-\varepsilon-\sigma-i(t+\tau)} dt.
$$
The denominator, say $D(t)$, is $C^\infty$ and $\gg 1+|t|$. Neglecting the constant factor, we have the sum of
$$
\Big[G(t) \frac{|t|^{m(\alpha-1/2)}}{D(t)}\Big]_{-T}^{T'} \ll T^{m(\alpha-1/2)} T^{-1},
$$
and of
$$
\int_{-T}^{T'} G(F)\frac{d}{dt}\Big(\frac{|t|^{m(\alpha-1/2)}}{D(t)}\Big) dt,
$$
computation justified if $\frac{|t|^{m(\alpha-1/2)}}{D(t)}$ is $C^1$. Since $\frac{d}{dt}(|t|^{m(\alpha-1/2)})$ is equal to $(Cst)|t|^{m/2-1+m\varepsilon}$ for $t$ positive or negative, this condition is satisfied if $m\ge 2$, which we will have to assume presently. Differentiating  the quotient $\frac{|t|^{m(\alpha-1/2)}}{D(t)}$, we find that the last integral is $\ll T^{m(\alpha-1/2)}T^{-1}$, as the first. The integral on the missing segment $[-T-\tau,-T]$ admits the same bound. Finally, the integral (2.19) is
$\ll (QT)^{m(\alpha-1/2)}y^{-1} T^{-1}$. Summing
 these contributions for $n>2(QT)^m/x$, we find that these terms contribute to (2.18) a term
\begin{equation}
O(x^{-\sigma}\ T^{-1}(QT)^{m/2}),
\end{equation}
already present in (2.17).

We now consider the remainder term
\begin{equation}
I_{-\varepsilon}^r = \varepsilon(\pi) x^{1-s} \sum_n \ b_n \ c(nx),
\end{equation}
where $nx\le 2(QT)^m$.

Again we follow \cite{FI}. We move the integration giving $c(y)$ from  $Re(w)=\alpha$ to  $Re(w)=\beta =\frac{1}{2}+\frac{1}{m}$. 

\textit{We assume} $m\ge 2$. (This is implicit in \cite{FI}.) However, we move only the segments with $1\le\vert t \vert\le T',T''$. We have to estimate the contributions of the vertical segment $Re(w)=\alpha$, $|t|\le 1$, and of the horizontal segments $Re(w) \in [\beta,\alpha]$, $t=T',-T''$ and $t=\pm1$. The vertical integral, $\gamma$ being the product of $D^{w-1/2}y^{-w}$ and of a bounded function (uniformly if the $c_i$ are fixed) is an $O(D^{1/2+\varepsilon} y^{-\alpha})$. The horizontal integrals for $t=T',-T''$ are dominated by $\int_\beta^\alpha(QT)^{m(u-1/2)} y^{-u}\ T^{-1}du$, dominated by
\begin{equation}
((QT)^{m(\alpha-1/2)}y^{-\alpha} + (QT)^{m(\beta-1/2)} y^{-\beta})T^{-1}
\end{equation}
by a computation similar to the one preceding (2.12). Since $y\le 2(QT)^m$, the first term is dominant. The horizontal integrals for $t=\pm1$ are dominated by $y^{-\beta}$. (It may occur that  $1-\sigma \in [\beta, \alpha] $ and $ \tau = \pm1$; in this case use a small indentation of the horizontal contour.) This is dominated by $(yT)^{-1} (QT)^{\frac{m}{2}}$: this amounts to $y^{1-\beta} \ll T^{-1} (QT)^\frac{m}{2}$, which follows from $y \ll (QT)^m$, $1-\beta= 1/2-1/m$ and $Q \gg 1$.

If $m\ge 2$,
$$
D^{1/2+\varepsilon}\asymp Q^{m(\frac{1}{2}+\varepsilon)} \ll (QT)^{m(\frac{1}{2}+\varepsilon)} T^{-1}.
$$

We now assume (as will be fulfilled later)

\vspace{2mm}

{\myhypo \textit{(Auxiliary)}. 
$$
(QT)^m \le x^2.
$$}

Since $y\ge1$,
$$\begin{array}{rl}
T^{-1}(QT)^{m(\frac{1}{2}+\epsilon)} y^{-\alpha} &\le (y T)^{-1} (QT)^{\frac{m}{2}} (QT)^{m\varepsilon}\\
\noalign{\vskip1mm}
&\le (y T)^{-1}(QT)^{m/2}x^{2\varepsilon}.
\end{array}
$$
Writing 
 $\overset{\beta+iT'}{\underset{\ \beta-iT''}{\displaystyle\intt}} (T',T''>1)$ for the integral with the segment $[-1,1]$ deleted, we have for $y\le 2(QT)^m$ :
\begin{equation}
c(y) = \frac{1}{2i\pi} \overset{\ \beta+iT'}{\underset{\beta-iT''}{\displaystyle\intt}}  \gamma(w) \frac{y^{-w}}{1-w-s} dw +O((y T)^{-1} (QT)^{m/2}x^{2\varepsilon}),
\end{equation}
compare \cite[(2.8)]{FI}.

However, the translation of integrals is justified only if we do not encounter zeroes of $1-w-s$ for $w=u+it$
in the upper rectangle $\beta \leq u \leq \alpha$,  $t \in [1;T']$ and the lower one. Since $w=1-s$, this means that $\sigma \in  ]0, 1/2-1/m]$; for the upper rectangle, $ t=-\tau$ (where we have chosen $\tau \geq0$) is impossible. For the lower rectangle, we get the value $t=\tau.$ 

In the latter case, we obtain a residue of order $\vert y^{-w} \gamma \vert \ll y^{-u} (Q\tau)^{m(u-1/2)}$ by Lemma 2.1. Since we will not pursue the dependence on $\tau,$ this is dominated by the remainder term in (2.24).  Indeed this amounts to

$$
y^{1-u} Q^{m(u-1/2)} \ll (QT)^{m/2} T^{-1} x^{2\varepsilon};
$$

Since $y \ll (QT)^m$,  this is implied by
 $$
 (QT)^{m(1-u)} Q^{m(u-1/2)} \ll (QT)^{m/2} T^{-1} x^{2\epsilon},
 $$
  so by 
  $$
  Q^{m/2}T^{m(1-u)} \ll T^{-1}(QT)^{m/2}
  $$
  i.e. \ $T^{m(1-u}) \ll T^{m/2-1}$, which is true since $1-u \leq 1/2-1/m$. (If $\sigma=1/2-1/m$, use a small indentation of the contour.)






\vskip2mm

Consider now the upper integral in (2.24). This is equal to
$$
\begin{array}{rl}
c_1(y) &= 	\dfrac{1}{2\pi}\ \omega_1\ Q\ y^{-\beta} \displaystyle \int_1^{T'} \dfrac{t(Qt/e)^{imt}y^{-it}}{1-s-(\beta +it)} (1+O(t^{-1}))dt\\
\noalign{\vskip 2mm}
&=\dfrac{i}{2\pi }\ \omega_1\ Q\ y^{-\beta} \displaystyle \int_1^{T'} (Qt/e)^{imt}y^{-it}\ dt +O((y T)^{-1}(QT)^{\frac{m}{2}}\cdot \log\ T)
\end{array}
$$
$(T'\ge 2)$, using that $y\ll (QT)^m$ and
\begin{equation}
\frac{1}{1-s-(\beta+it)} = \frac{i}{t}(1+O(t^{-1})).
\end{equation}
(We have used that $T' \asymp T$, $T'' \asymp T$).

Since $T\le x$, the remainder term is dominated by the last remainder term, $(y T)^{-1}(QT)^{m/2} x^{2\varepsilon}$. The lower integral yields

$$
\begin{array}{rl}
c_2(y) &= 	\dfrac{1}{2\pi}\ \omega_2\ Q\ y^{-\beta} \displaystyle \int_1^{T''} \dfrac{t(Qt/e)^{-imt}y^{it}}{1-s-(\beta -it)} (1+O(t^{-1}))dt\\
\noalign{\vskip 2mm}
&=\dfrac{-i}{2\pi }\ \omega_2\ Q\ y^{-\beta} \displaystyle \int_1^{T''} (Qt/e)^{-imt}y^{it}\ dt +O((y T)^{-1}(QT)^{\frac{m}{2}}\cdot \log\ T).
\end{array}
$$

We are now able to reduce the computation to that in \cite{FI} for $c(y)$. To be precise, in the case where $\omega_2 =\bar{\omega_1}=\bar{\omega}$,
 $c(y)$ is given by \cite[(2.8)]{FI} :
$$
c(y) = \pi^{-1}\ y^{-\beta}\ Q\ Re(i\omega \int_1^T (QT)/e)^{int} y^{-it}dt)+ \mathrm{remainder\ term.}
$$
(The parentheses are missing in the original.) Here this is replaced by
$$
\begin{array}{rl}
c(y) &= c'_1y)+c'_2(y) + \ \mathrm{remainder\ term},\\
c'_1(y) &=\dfrac{1}{2\pi} y^{-\beta}\ Q\ i\omega_1
\ \displaystyle\int_1^{T'}(Qt/e)^{imt}\ y^{-it}\ dt,\\
c'_2(y) &=\dfrac{1}{2\pi} y^{-\beta}\ Q\ (-i\omega_2)
\ \displaystyle\int_1^{T''}(Qt/e)^{-imt}\ y^{-t}\ dt
\end{array}
$$
We majorize the contributions of the intervals $[T',T]$ and $[T,T'']$ : write $C'_1(y),C'_2(y)$ for the integrals with end value $T$.

\vspace{2mm}

{\monlem 

\vspace{2mm}

$(1)$ $|c'_1(y) +c'_2(y)-C'_1(y)-C'_2(y)| \ll y^{-\beta} \ Q\ \tau$

\vspace{2mm}

$(2)$ $\displaystyle\sum_{nx\le 2(Q\Gamma)^m} b_n\ \tau\ Q\ (nx)^{-\beta} \ll \tau (x T)^{-1} (QT)^{m/2}x^\varepsilon$.
}

\vspace{2mm}

The first inequality is clear. The second follows by an easy summation, using that $b_n\ll n^{\varepsilon'}$ for any $\varepsilon'>0$.

Now the integrals $C'_1(y)$, $C'_2(y)$, and their sum over $n\le \frac{2(Q\Gamma)^m}{x}$, are estimated in \cite{FI}. For $1\le T\le x$ let
$$
B_0(x,T) = \sum_{nx\le(QT)^m} b_n\ n^{-\frac{m+1}{2m}}\ \Big\{\omega_1 \ e^{3i\pi/4-\frac{im}{Q}(nx)^{1/m}} + \omega_2\ e^{-\frac{3i\pi}{4}+\frac{im}{Q}(nx)^{1/m}}\Big\}.
$$
We have transformed the expression of $B_0(x,T)$ in \cite{FI} according to $(\omega,\bar{\omega}) \rg (\omega_1,\omega_2)$. Let
\begin{equation}
J_{-\varepsilon}^r(x) = \varepsilon(\pi)\ x^{1-s} \sum_n b_n\ c_{FI}(nx),
\end{equation}
where, as before Lemma 2.2, we have replaced $c'_1$, $c'_2$ by $C'_1$, $C'_2$; thus $c_{FI}$ is the function "c" of Friedlander-Iwaniec (with $\omega, \bar{\omega}$ replaced by $\omega_1, \omega_2$.) Note the factor $x^{1-s}$ in $(2.26)$, as opposed to $x$ in \cite{FI}.

We now have, using the estimate in \cite{FI} for $C'_1, C'_2$ :

\vspace{2mm}

{\monlem \rm \cite{FI} .- For $1\le T\le x$, with $(QT)^m \le x^2$, 
$$
\begin{array}{l}
J_{-\varepsilon}^r (x) = \varepsilon(\pi) \ x^{-s}\Big(\dfrac{2Q}{\pi m}\Big)^{1/2} x^{\frac{m-1}{2m}}B_0(x,T)\\
\noalign{\vskip2mm}
\hspace{1cm} +O(x^{1-\sigma} \ T^{-1/2}(QT)^{-m/2}) +O(x^{2\varepsilon-\sigma}\ T^{-1}(QT)^{m/2}).
\end{array}
$$}
\vspace{2mm}

We still have to replace $J_{-\varepsilon}^r$ by $I_{-\varepsilon}^r$. The difference, according to\break Lemma~2.2, is dominated by the second remainder term in Lemma 2.3 (multiplied by $\tau$, but we assume here $s$ fixed.). We therefore obtain the same estimate.

Collecting the remainder term (2.17), and noting that
$$
x\ T^{-1/2}(QT)^{-m/2} \ll x\ Q^{-m/2}\ T^{-\frac{1}{2}-\frac{m}{2}} \ll x\ T^{-1}
$$
since $Q\ge \frac{1}{2\pi}$, we obtain :
\begin{eqnarray}
{A_s(x)\! = \!\frac{\kappa}{1\!-\!s} \ x^{1-s}\! +\! L(s)\! +\! \varepsilon(\pi)x^{-s} \big(\frac{2Q}{\pi m}\big)^{1/2}\, x^{\frac{m-1}{2m}} \cdot}\\
\cdot B_0(x,T)\! +\! R(x,T),\nonumber
\end{eqnarray}
\begin{eqnarray}
R(x,T) &\ll& x^{-\sigma} \{x^{1+\varepsilon} +(QT)^{\frac{m}{2}} x^{2\varepsilon}\}T^{-1}\\
\noalign{\vskip2mm}
&\ll& x^{-\sigma+\varepsilon} \{x+(QT)^{\frac{m}{2}} \ x^\varepsilon \} T^{-1}.
\end{eqnarray}
This is essentially the expression in \cite[p. 499]{FI}, the remainder term being of course multiplied by $x^{-\sigma}$.

Friedlander and Iwaniec then set $N=x^{-1}(QT)^m$. Assuming $N\le x$, $R(x,T)\ll x^{1-\sigma+\varepsilon}T^{-1}$. We have $B_0(x,T)=B(x,N)$ where \cite[(1.12)]{FI}
\begin{eqnarray*}
B(x,N) = \displaystyle\sum_{n\le N} b_n n^{-\frac{m+1}{2m}} \{ \omega_1\ e^{-i\pi/4-2i\pi m(\frac{nx}{D}})^{1/m} +\\
\noalign{\vskip2mm}
\omega_2\ e^{i\pi/4+2i\pi m(\frac{nx}{D})^{1/m}}\}.
\end{eqnarray*}

Obviously the estimates for $B(x,N)$ and for its variant given in \cite{FI} by (1.12) are the same. However we must check that $1 \leq N\le x$, as in \cite[Theorem~1.2]{FI} are allowed values. For $1\le T\le x$,
$$
Q^m/x \le N \le Q^m\ x^{m-1} = D/(2\pi)^m\cdot x^{m-1}.
$$
Since $m\ge 2$, $N \leq x$ (essentially) satisfies the upper bound. However $N=1$ is allowed only if $x\ge D/(2\pi)^m$, not for $x\ge D^{1/2}$ as stated. Note also that Auxiliary hypothesis~2.2 is satisfied. Under this assumption, we get the statement of Theorem 1.2 in \cite{FI}, the remainder term being multiplied by $x^{-\sigma}$. Namely, for $x \geq D$,
\begin{eqnarray}
{A_s(x)\! = \!\frac{\kappa}{1\!-\!s} \ x^{1-s}\! +\! L(s)\! +\! \varepsilon(\pi)x^{-s} \big(\frac{Q}{2\pi m}\big)^{1/2}\, x^{\frac{m-1}{2m}} \cdot}\\
\cdot B(x,N)\!  + O(D^{1/m} N^{-1/m} x^{-\sigma+  \frac{m-1}{m}+\varepsilon})\nonumber
\end{eqnarray}

for  $1 \leq N \leq  (2\pi)^{-m }x$.

 We want to evaluate  this, in order to get Proposition~1.1 in \cite{FI}, at $N=D^{\frac{1}{m+1}}\ x^{\frac{m-1}{m+1}}$. We must therefore have
$$
D^{\frac{1}{m+1}}\ x^{\frac{m-1}{m+1}} \ge x^{-1}\ Q^m = x^{-1}\ D(1/2\pi)^m,
$$
which yields by a simple computation
\begin{equation}
x^2 \ge D(1/2\pi)^{m+1}.
\end{equation}
Thus the stated condition, $x^2\ge D$, is adequate. Under this assumption, Friedlander-Iwaniec show that the full remainder term for $s=0$:
$$
\varepsilon(\pi) \big(\frac{Q}{2\pi m}\big)^{1/2} x^{\frac{m-1}{2m}}
B(x,N)\!  + O(D^{1/m} N^{-1/m} x^{  \frac{m-1}{m}+\varepsilon}) 
$$

is $\ll D^{\frac{1}{m+1}} x^{\frac{m-1}{m+1}+\varepsilon}$.

 We have introduced another condition in the proof, i.e. $T\ge 2\tau$, so for $N=x^{-1}(QT)^m$ :
$$
N\ge x^{-1}(D/(2\pi)^m)(2\tau)^m
$$
which yields
\begin{equation}
x\ge D^{1/2}\Big|\frac{\tau}{\pi}\Big|^{\frac{m+1}{2}}.
\end{equation}

\vspace{2mm}

{\montheo Assume $s=\sigma+i\tau$, $0<\sigma<1$. Then, if
$$
\begin{array}{c}
x^2 \ge \Max(D(1/{2\pi})^{m+1}, D \vert(\tau / \pi) \vert^{m+1}),\\
\noalign{\vskip2mm}
A_s(x) = \dfrac{\kappa}{1-s}\ x^{1-s} + L(s) +O\big(D^{\frac{1}{m+1}}\ x^{\frac{m-1}{m+1}-\sigma+\varepsilon}\big)
\end{array}
$$
the implicit constant depending only on $\varepsilon$, $\pi_\infty$ and $\tau$.
}

\vspace{2mm}

\noindent \textbf{{Corollary.}} \textit{Under the same assumptions,}
$$ 
H_s(x)=O\big(D^{1/m+1}\ x^{-\frac{2}{m+1}+\varepsilon}\big)
$$
\textit{if} $L(s)=0$.

\vspace{2mm}

As expected, we have improved on the estimates on $H_s$ (before Lem\-ma~1.1)~: the new estimates do not depend on~$\sigma$.

\section{Automorphic $L$--functions ~ (mostly) ~ of degree  $\le2$ over~$\Q$}

\subsection{} In this chapter we use the results of Chapter~2 to extend Theorem~1.2. We consider $L$--series associated to automorphic representations of $GL(m,\A_\Q)$ $(m=1,2)$ ; moreover we will consider the integrals of $\big|\frac{L(w)}{w-s}\big|^2$ on vertical lines $Re(w)=\tau$ where $\tau\not= 1/2.$\footnote{Thus $\tau$ no longer denotes $Im(s)$ as in Chapter~2 !}

In \S~3.4, we explore the consequences of Theorem~2.2 for the convergence of the series $L(s,\pi)$ associated to a cuspidal representation of $GL(m)$. In \S~3.5, we discuss this problem, and our main conjecture, in light of standard conjectures about $L$--functions. Finally, in \S~3.6, we briefly inquire whether there exist \textit{upper bounds}, as well as lower bounds, for our quadratic integrals.

\subsection{} Consider first the case of $\zeta_\Q$, or of a classical Dirichlet series $L(s,\chi)$. (We will always consider \textit{primitive} characters.) In this case we saw that $\Mm H_s$ (if $\zeta(s)$ or  $L(s, \chi)=0$) is defined (and holomorphic) for $0<Re(w)<1$.\footnote{Strictly speaking we did not do this for $L(s,\chi)$ but the same estimate $H_s(x)\ll x^{-1}$ is obtained by extending the computation in \S~1.3.} We can apply the formula of Mellin--Parseval in this domain.

Assume $f(x)$ $(x\ge 0)$ is a function bounded in $0$ and such that $f(x) \ll x^{-\lambda+\varepsilon}$, $\lambda>0$ (for any $\varepsilon>0$) if $x\rg \infty$. Then $\Mm f(w)$ is defined and holomorphic for $0<Re(w)<\lambda$. Moreover, for $\tau<\lambda$,
$$
\int_0^\infty\ | x^\tau \ f(x)|^2\ \frac{dx}{x} < \infty.
$$
We are then in the domain of Titchmarsh's Theorem 71 \cite{T2} : $\Mm f$ is $L^2$ on the line $Re(w)=\tau<\lambda$, and
$$
\int_0^\infty\ |x^\tau\ f(x)|^2 \frac{dx}{x} = \frac{1}{2\pi} \int_{-\infty}^{+\infty} \ |\Mm f(\tau+it)|^2dt.
$$

\vspace{2mm} 

 For $-2\le a\le 2$, 

$$
\int_1^2\ x^{a-1}\ dx \ge \frac{3}{8}.
$$


Using the argument for the proof of Theorem~1.2 $(i)$, we deduce :

\vspace{2mm} 

{\maprop Assume $s$ is a zero of $L(s)=\zeta_\Q(s)$ or $L(s,\chi)$. Then, for any $\tau\in ]0,1[$,
$$
\int_{-\infty}^{+\infty} \Big|\frac{L(\tau+it)}{\tau+it-s}\Big|^2\ dt > \frac{3\pi}{4}\  - \frac{4\pi}{|1-s|}
$$
(the remainder term being absent for a non--trivial Dirichlet character).}

\vspace{2mm} 

Consider now an automorphic representation $\pi$ of $GL(2,\A_\Q)$ verifying the assumptions of \S~2.1. (In particular $\pi_\infty$ is self--dual). This includes the case of $L_F(s,\chi)$ where $F/\Q$ is quadratic and, if $F$ is  real, $\chi$ is an Artin character and, if $F$ is imaginary, $\chi$ is an algebraic Hecke character (possibly multiplied by $|\ |_F^{it}$), and in particular the case of $\zeta_F$. By the Corollary to Theorem~2.1, $H_s(x) \ll x^{-\frac{2}{3}+\varepsilon}$. In particular (see after Proposition~1.2) we see that $\Mm H_s(w)$ is convergent for $\tau=1/2$, and that $H_s \in L^2$ (Proposition~1.2 $(i)$) follows from the growth estimate.

In fact $\Mm H_s(w)=\frac{L(1-w)}{1-w-s}$ is convergent in the domain $\tau<\frac{2}{3}$, and is $L^2$ on the vertical lines. Arguing as before, we have~:

\vspace{2mm} 

{\maprop Assume $\pi$ is a tempered, automorphic representation of $GL(2,\A_\Q)$ with $\pi_\infty$ self--dual. Let $s$ be a zero of $L(s,\pi)$, with $0<\sigma<1$. 
Then, for any $\tau>\frac{1}{3}$, $\frac{L(\tau+it)}{\tau+it-s}$ is $L^2$ and
$$
\int_{-\infty}^{+\infty} \Big| \frac{L(\tau+it)}{\tau+it-s}\Big|^2 dt > \frac{3\pi}{4} - 4\pi \Big|\frac{\kappa}{1-s}\Big|.
$$
(For $\tau=1/2$, the constant $\frac{3\pi}{4}$ can be replaced by $2\pi \log 2$.)
Moreover, for $\sigma=1/2 $:
$$
\int_{-\infty}^{+\infty} \Big| \frac{L(1/2+it)}{1/2+it-s}\Big|^2 dt > \pi \rm \  log 2.
$$}

\vspace{2mm} 

The argument for the second part is as follows. We seek a lower bound on

$$
\int_0^2 \vert H_s(x) \vert^2 dx = \vert \frac{\kappa}{1-s} \vert^2 +  \int_1^2 \vert x^{s-1}-\frac{\kappa}{1-s} \vert^2 dx 
$$
Let $f(x)=x^{s-1}, g(x) =\frac{\kappa}{1-s} $ on $[1,2]$. The full integral is then equal to $\vert\vert  f-g\vert\vert^2 + \vert\vert g \vert\vert^2$ ($L^2$-norm on [1,2].) However this is greater than $\vert\vert f \vert\vert^2-2 \vert\vert f\vert\vert. \vert\vert g \vert\vert +2 \vert\vert g \vert\vert^2 \geq \frac{1}{2} \vert\vert f \vert\vert^2 = 1/2 \ \rm log \ 2$, whence the result.

It is interesting to compare this with  the subconvex estimate. Let $\mu(\tau)$ be defined by
$$
|L(\tau+it| \ll t^{\mu(\tau)+\varepsilon}
$$
and minimal. Then $\mu(\frac{1}{2}) \le \frac{1}{2}-\delta$ $(\delta>0)$, $\mu(0) =1$ and therefore $\mu(\tau) \le 1-(1+2\delta)\tau$ for $\tau\in ]0,\frac{1}{2}]$. In particular $\frac{L(\tau+it)}{t}$ $(t\ge 1)$ is $L^2$, as a consequence of the convexity estimate, if $\tau>\frac{1}{2(1+2\delta)}$. Since $\delta$ is small, this does not imply the convergence of the quadratic integral in the range $\tau >\frac{1}{3}$.

\subsection{The case $s=0$\protect\footnote{This was the case considered by Titchmarsh \cite[p. 321]{T1}. 
Titchmarsh, studying the divisor problem, considers $L(s)=\zeta_\Q(s)^2$, a case we have excluded here. He uses the Mellin transform in the same fashion.}}

In this section we prove Theorem C and Theorem D in the context of this paper, i.e. when the integrals are convergent. We make the assumptions in \S~2. Thus
$$
A_0(x) = \sum_{n\le x} a_n = \kappa \, x +O\big(D^{\frac{1}{m+1}} \, x^{1-\frac{2}{m+1}+\varepsilon}\big),
$$
a priori for $x\ge D^{1/2}$. (However $|a_n|\le \tau_m(n)$ and so
$$
\sum_{n\le x} |a_n| \ll x(\log x)^{m-1},
$$
cf. \cite[p. 23]{FI}, and a trivial computation shows that, if $\kappa=0$, this estimate is true for all $x \geq 1$.)

We have $H_0(x) = x^{-1}A_0(x)-\kappa$; note that $(s=0)$ is not necessarily a zero of $L(s)$. As before we see that
$$
\Mm H_0(w) = \int_1^\infty H_0(x) \, x^{w-1} dx
$$
is convergent if $\tau={Re}(w) < \frac{2}{m+1}$. We will now assume that $m=1,2$ and consider $\tau=1/2$. 
Assume first that $\kappa=0$. We can proceed as in \S~1.3 and compute directly $\Mm H_0$ :
$$
\begin{array}{l}
\displaystyle\int_1^T x^{-1} \displaystyle\sum_{n\le x} a_n \, x^{w-1}dx\\
= \displaystyle \sum_{n\le T}a_n \int_n^T x^{w-2} dx\\
= \dfrac{1}{w-1} \Big\{ T^{w-1} \displaystyle\sum_{n\le T} a_n - \sum_{n\le T} a_n\, n^{w-1}\Big\}.
\end{array}
$$

The first bracketed term is $\ll T^{\tau-\frac{2}{m+1}}\rg 0$ $(T\rg \infty)$, and the second is $\sum\limits_{n\le T}a_n \, n^{w-1}\rg L(1-w)$ by Theorem~2.2 (cf. also Theorem~3.1). Thus
$$
\Mm H_0(w) = \frac{L(1-w)}{1-w}
$$
for $\tau<\frac{2}{m+1}$ and in particular for $\tau=\frac{1}{2}$.

In the case where $L(s)$ has a pole (at $s=1$) we must proceed as in \S~1.4. We have
$$
A_0^*(x) =\frac{1}{2i\pi} \int_{c-i\infty}^{c+i\infty} L(w) \frac{x^w}{w}dw
$$
for $c>1$,
$$
x^{-1}A_0^*(x) = \frac{1}{2i\pi} \int_{c'} L(1-w) \frac{x^{-w}}{1-w}dw
$$
for $c'=1-c<0$, by a change of variables. We move the integral to $c''>0$, $c''<1$ (e.g. $c''=\frac{1}{2}$) and deduce
$$
x^{-1} A_0^*(x) = \frac{1}{2i\pi} \int_{1/2} \frac{L(1-w)x^{-w}}{1-w} dw+\kappa
$$
using the same estimates as in \S~1.4, taking account of the residue at $(w=0)$. Thus
$$
H_s^*(x) = x^{-1} A_0^*(x)-\kappa
$$
is given by the vertical integral.

The first computation yields Theorem C by the Parseval formula, since
$$
\begin{array}{l}
\displaystyle\int_1^\infty |x^{1/2} H_0(x)|^2 \frac{dx}{x} = \int_1^\infty |H_0(x)|^2 dx\\
\displaystyle\ge \int_1^2 x^{-2}dx = \frac{1}{2}.
\end{array}
$$
The second case yields
$$
\frac{1}{2\pi} \int \Big|\frac{L(1/2+it)}{1/2+it}\Big|^2 dt \ge\int_0^1 |\kappa|^2 dx + \int_0^2 |x^{-1}-\kappa|^2 dx 
$$

The argument given in the proof of Proposition 3.2 shows that this sum of integrals is larger than $1/4$,
proving Theorem D.

\vspace{2mm}

Note the special case of Theorem D :

\vspace{2mm}

{\montheo Assume $F$ is a quadratic field. Then \break $\int \Big|\dfrac{\zeta_F(\frac{1}{2}+it)}{1/2+it}\Big|^2 dt > \pi /2$.
}

\vspace{2mm}

We now consider the case of $m=1$, and of a non--trivial Dirichlet character $\chi$. We assume that $\chi$ is a character $\mod q$, $q$ being prime. Here the Friedlander--Iwaniec estimate is simply replaced by the Polya--Vinogradov inequality (Schur's version) $A_0(x) \leq q^{1/2} $. The vertical integral is bounded below by $\pi$. However, under the Lindel\"of hypothesis \cite[Cor.~5.20]{IK} the integrand admits a majoration by $\dfrac{(q|s|)^{2\varepsilon}}{1/4+t^2}$, so the integral $\ll q^{2\varepsilon}$. We now consider the quadratic character $\chi$ of $(\Z/q)^\times$.

Let $\nu=\nu(q)$ be the least quadratic non--residue. We have
$$
\begin{array}{rl}
\displaystyle \int_1^\infty |H_0(x)|^2 dx &\ge  \int_1^\nu (x^{-1}[x])^2 dx\\
&\ge \displaystyle \int_1^\nu \Big(1-\dfrac{2}{x}\Big) dx =\nu -2 \log \nu
\end{array}
$$
using that $[x]\ge x-1$. Thus
$$
\nu-2 \log \nu \ll q^{2\varepsilon},
$$
i.e. $\nu\ll q^\varepsilon$. Thus this conjecture of Vinogradov follows directly from the \textit{Lindel\"of} hypothesis (including its dependence on $q$.) The proof of Ankeny \cite{ An}  depends on the generalised Riemann hypothesis.

\subsection{} We now consider the convergence abscissa of a cuspidal representation $\pi$ of $GL(m,\A)$. We  assume $\pi$ tempered and $\pi_\infty$ self--dual.

\vspace{2mm} 

{\montheo 

$(i)$ The abscissa of absolute convergence $\sigma_a(\Pi)$ of the Euler product for $L(s,\pi)$ satisfies $0\leq \sigma_a(\Pi) \leq 1$ .

$(ii)$ The abscissa of absolute convergence of the Dirichlet series for $L(s,\pi)$ satisfies $ 0 \leq\sigma_a \leq 1$.

$(iii)$ The abscissa of convergence $\sigma_c$ of the $L$--series $L(s,\pi)$  verifies $ \sigma_c \le 1-\frac{2}{m+1}$.}

\vspace{2mm}

That $\sigma_a(\Pi) \leq 1$  is due to Jacquet--Shalika \cite{JS}. Let $A_p$ be the Hecke matrix of $\pi$ at an unramified prime $p$, so $a_p=\textrm{ trace}(A_p)$. If  \ $\sum \vert a_p\vert p^{-\sigma} < \infty$ for $\sigma<0$, $\vert \textrm{ trace}(A_p)\vert^2 p^{-2\sigma} \leq C$, which implies that $\sum \vert \textrm{ trace}(A_p)\vert^2 p^{-1} < \infty.$  This contradicts the divergence of the logarithm of the Rankin L-function at $s=1$ \cite{JS}. Thus $\sigma_a(\Pi) \geq 0$.  Again, by \cite{JS}, $\sigma_a\le 1$. Then $\sigma_a\ge 0$ follows  from the usual argument : if $L(\sigma,\pi)$ were absolutely convergent for $\sigma<0$, we would have

$$
\sum_p \ \sum_{\alpha\ge 1}
 \ |a_{p^\alpha}|\ p^{-\alpha\sigma} < \infty
 $$
 and this would imply that the Euler product is convergent. 
 
 Finally, Theorem 2.2 implies that $\sum\limits_{n\le x} a_n\ n^{-s}$ converges to $L(s)$ if $\sigma<\frac{m-1}{m+1}$.
 
\textbf{Remark}.- If $\pi$ verifies the generalised Sato-Tate conjecture, it is easy to see that $\sigma_a(\Pi)=1$. This applies, in particular, when $m=2$ and $\pi$ is associated to classical modular forms.  
 
 Anticipating on the next paragraph, we note that a stronger result should hold under a Conjecture of Friedlander--Iwaniec \cite[Conjecture 2]{FI}, which would imply (we forget the dependence on $D$)~:
 \begin{equation}
 A_0(x) \ll x^{\frac{m-1}{2m}+\varepsilon}.
 \end{equation}

One should expect that the proof in Chapter~2 would then extend to yield
\begin{equation}
A_s(x) - L(s) \ll x^{\frac{m-1}{2m}-\sigma+\varepsilon},
\end{equation}
whence
\begin{equation}
 \sigma_c \le \frac{1}{2}-\frac{1}{2m}.
\end{equation}

Assume $\sigma_c\le\frac{1}{2}$. By a classical result $\mu(\sigma_c)\le 1$ \cite[9.33]{T3}. On the other hand \cite[9.41]{T3} , $\mu(\sigma)=0$ for $\sigma>1$ and the functional equation implies $\mu(\sigma)= m(1/2-\sigma)$ for $\sigma<0$ and therefore  $\mu(\sigma)\geq m(1/2-\sigma)$ by convexity for $\sigma \leq 1/2.$ Thus we see (unconditionally) that
\begin{equation}
\sigma
_c \ge \frac{1}{2}-\frac{1}{m}.
\end{equation}

We do not know, even hypothetically, what should be the correct value of~$\sigma_c$.

\subsection{} We now describe, in any degree $m$, what should be expected of the square--integrability  of $H_s$ and $\frac{L(1-w)}{1-s-w}$, in which case the lower bound would be effective. It is interesting to consider, as in \S~3.2, values of $Re(w) \not=\frac{1}{2}$. As we will see, we need the full force of the Friedlander--Iwaniec conjecture, i.e.~(3.1). [Here $L(s)$ is an arbitrary $L$--function  $L(s,\pi)$ of degree $m$ over $\Q$, where $\pi$ satisfies the conditions in~\S~2.1.]

Consider first the convergence of $\Mm H_s(w)$. One expects (3.2), which implies
$$
H_s(x) \ll x^{-1/2-\frac{1}{2m}+\varepsilon}.
$$
Thus $\Mm H_s(w)$ is well--defined for $\tau=Re(w) <\frac{1}{2}+\frac{1}{2m}$. Moreover, $\Mm H_s(w)$ will be $L^2$ on the line $Re(w)=\tau$ if (see \S~3.2)
$$
\int_0^\infty \ x^{2\tau}\ |H_s(x)|^2\ \frac{dx}{x} < \infty,
$$
which gives the same condition. Thus we would have by the substitution $w\rg 1-w$~:

\vspace{2mm}

\textit{Assume $\tau>\frac{1}{2}-\frac{1}{2m}$. Then, for any zero of $L(s)$, $\Mm H_s(1-w)$ is defined, and equal to $\frac{L(w)}{w-s}$ for $Re(w)=\tau$. Moreover,}
$$
\frac{1}{2\pi} \int \ \Big| \frac{L(\tau+it)}{\tau+it-s} \Big|^2 dt = \int_0^\infty\ x^{2-2\tau} |H_s(x)|^2\ \frac{dx}{x}.
$$

\textit{Both integrals are finite. In particular, $\frac{L(\tau+it)}{1+|t|}$ is~$L^2$.}

Assume $\frac{1}{2}-\frac{1}{2m}<\tau \le \frac{1}{2}$. Then, under the generalized Lindel\"of conjecture,
$$
\mu(\tau) = m(\frac{1}{2}-\tau) <\frac{1}{2}.
$$
This implies that $\frac{L(\tau+it)}{\tau+it-s}\in L^2$, in conformity with the previous (conjectural) description.

It is easy to see that the estimate coming from the generalised Riemann hypothesis, namely, cf. \cite[p. 500]{FI} ; \cite[p. 116]{IK}:
\begin{equation}
A_0(x) = \kappa\ x+O(x^{1/2+\varepsilon})
\end{equation}
does not allow, in general, for reaching the line $\tau=\frac{1}{2}$.

\vskip2mm

\textbf{Remark}.- According to Iwaniec-Kowalski \cite[p.86]{IK}, the effective estimate (of Landau or Friedlander-Iwaniec) for $A_0$ can be improved. Any small improvement would allow one to prove the results in this paper (with finite integrals) for $m=3$.

\subsection{}

    Since we expect an absolute lower bound, when $L(s)=0$, on the quadratic integral of $\frac {L(w)}{w-s}$ on the line $Re(w)= 1/2$, it is now a natural question whether this integral will tend to $\infty$ when $s \rightarrow \infty$ in the critical strip. Even in the case of the Riemann zeta function, this seems a difficult problem. We will only make a few remarks. As in chapter 1, let
    
    \vskip 2mm
    
    $$
    I(s) =\int_{-\infty}^{+\infty} \Big|\frac{\zeta(\frac{1}{2}+it)}{1/2+it-s}\Big|^2 dt.
    $$
    
    Suppose that $I(s)$ is bounded. It is trivial that $ \int_a^b \frac{\zeta(\frac{1}{2}+it)}{1/2+it-s} dt \rightarrow 0$ for $s\rightarrow \infty$ for any finite interval $[a,b]$. Since the linear combinations of characteristic functions of such intervals are dense, we see that under this assumption
    
    \vskip 2mm
    $$
     \frac{\zeta(\frac{1}{2}+it)}{1/2+it-s} ~~ and ~~~ H_s ~~~ tend~~ to ~~0  ~~weakly ~~ in~~ L^2.
   $$
   
   \vskip 2mm
   
      The boundedness of $I(s)$ would also have a puzzling consequence. We consider only values of $s$ on the critical line. Writing $z(t)$ for $\zeta(1/2+it)$ and $s=1/2+i\tau$,  we would then have
      $$
      \int_{-\infty}^{+\infty} \Big|\frac{z(t)}{t-\tau}\Big|^2 dt \leq C.
     $$
     
     We consider the integral on a small interval around $\tau$. We have $z(\tau+h) = h z'(\tau) +\frac{h^2}{2} (z_{re}''(\tau+\theta_1 h) + iz_{im}''(\tau+\theta_2 h))$ where we have decomposed $z''$ in its real and imaginary parts.
      If $h\ll \tau$ and $\mu$ is an order of $\zeta$ on ($\sigma=1/2$), for example $\mu= 1/6$, $z''$, which is essentially $\zeta''$, is dominated by $\tau^{\mu+\varepsilon}$ by the residue formula. Replacing $\varepsilon$ by $2\varepsilon$ and taking $\tau$ sufficiently large, we deduce that for $u \leq h$: 
     $$ 
   \vert \frac{ z(\tau+u)}{u} \vert  \geq \vert z'(\tau)\vert  - \frac{u}{2} \tau^{\mu+\varepsilon},
    $$
    whence
    $$
      \vert \frac{ z(\tau+u)}{u} \vert^2 \geq \vert z'(\tau)\vert^2 - u \tau^{\mu+ \varepsilon} \vert z'(\tau) \vert
     $$
     and
     $$
         \int_{0}^{h} \Big|\frac{z(\tau+u)}{u}\Big|^2 du \geq \vert z'(\tau) \vert(h\vert z'(\tau) \vert - \frac{ h^2}{2} \tau^{\mu+\varepsilon}).
     $$
     
      Choosing $h= \frac{ \vert z'(\tau) \vert}{\tau^{\mu+\varepsilon}}$, so clearly $h \ll \tau$,  we see that the integral  is larger than
        $\frac{1}{2} \vert z'(\tau) \vert^3  \tau^{-\mu-\varepsilon}.$ Since  this is uniformly bounded, we deduce 
        
   {\maprop Assume I(s) is uniformly bounded for s a zero on the critical line. Then, for any such zero $s=1/2+i\tau$, 
$$
\zeta'(s) \ll \tau^{\mu/3+\varepsilon}.
$$} 
Of course the zeroes are too sparse to use this argument to improve the estimate on $\zeta (s)$ on the critical line. Although the estimate on $\zeta'$ is undoubtedly true, it seems unlikely that it can be so inferred.

\eject

\eject

\begin{center}
\textbf{\large Appendix. Variances for $\lambda_\pi(n)$ and $L(s,\pi)$}\\

\vskip3mm

\textbf{\large L. Clozel and P. Sarnak}
\end{center}


\section*{Introduction }
 
 The universal lower bounds for variances of $L$--functions on their critical lines that are established in this paper are based on a Parseval identity. The $L^2$--duality allows one to express a variance on the one side in terms of a quadratic form on the other side, facilitating  an estimation. In this appendix we examine variances for the variations of $\lambda_\pi(n)$ in $n$ and $L(\frac{1}{2}+it,\pi)$ in $t$, and place the results of this paper in this more general context.
 
 We restrict to $\pi$'s which are cuspidal on $GL_m/\Q$, $m\ge 1$, including the case that $\pi$ is the trivial representation and $L(s,\pi)=\zeta(s)$. We also specialize the point $s$ to be $1/2$ in the notation of the paper, which we adopt. The analytic conductor $c(\pi)$ (see \cite{IS}) is defined to be
  $$
 D\cdot \prod_{j=1}^m (1+|c_j|)
 \leqno(1)
 $$ 
 where $D=D(\pi)\in \N$ is the usual conductor of $\pi$ (see \S 2.1) \footnote{The unspecified references are to the main paper.} and $c_j=c_j(\pi)$ are as in (2.2). So $D(\pi)$ measures the finite ramification of $\pi$ and the $c_j$'s correspond to  the Archimedian component $\pi_\infty$ of $\pi$. $c(\pi)$ is a measure of the complexity of $\pi$ and it enters when estimating quantities associated with $\pi$. Note that for $m$ fixed, the number of $\pi$'s with $c(\pi)$ at most a given limit, is finite. In fact this count satisfies a ``Weyl--Schanuel'' asymptotic law (see \cite{BM}).
 
 In the main paper, only the arithmetic conductor occurs. We note however that the Friedlander-Iwaniec estimate (Prop. 1.1) can likely be formulated in terms of the analytic conductor, and doing so would be worthwhile.
 
 Let $L(s,\pi)= \sum_{n=1}^{\infty} \lambda_{\pi}(n)n^{-s}$. We are interested in the fluctuations of the function $\lambda_{\pi}(n)$, or rather of the summatory function $\sum_{n \leq x} \lambda_{\pi}(n)$ ($x\geq 1$.)

On the dual side are the fluctuations of the function $L(\frac{1}{2}+it,\pi)$ as functions of $t\in\R$. We fixate on their fluctuations about their value at $t=0$. In Theorems A and B, this corresponds to $s=1/2$; however we do not want to assume that $L( \frac{1}{2},\pi)=0$. The variance $V$ introduced in this paper is defined by :
$$
\begin{array}{lr}
 \hspace{1cm} &V(\pi) = V(\frac{1}{2},\pi) := \displaystyle\int_{-\infty}^\infty \frac{|L(\frac{1}{2}+it,\pi)-L(\frac{1}{2},\pi)|^2}{t^2} dt.
\end{array}\leqno(2)
$$
As discussed in detail in the paper, $V(\pi)$ is expected to be finite. We will allow it to be infinite; as a consequence, we make no use of any unproven hypotheses.

The key Parseval relation proved in the paper may be extended to prove (assuming $L(s,\pi)$ has no pole) that :
$$
V(\pi) = 2\pi \int_0^\infty \Big| \sum_{n\le x} \frac{\lambda\pi(n)}{\sqrt{n}} -L(\frac{1}{2},\pi)\Big|^2 \frac{dx}{x}.
\leqno(3)
$$
The meaning in (3) is that both sides are finite or infinite together, and if finite they are equal. Note that (3) expresses variance $V(\pi)$ in terms of mean--square of $\lambda_\pi(n)/\sqrt{n}$ with $n$ varying over $\N$.

The means and variance of our fluctuating variable give us a good picture of its distribution. How do these vary with $\pi$ ? The size of the central value $L(\frac{1}{2},\pi)$, is  a much studied problem (see \cite{Mu} for a survey). Any improvement in the exponent of the bound ($m$ fixed)
$$
L(\frac{1}{2},\pi) \ll c(\pi)^{1/4}
\leqno(4)
$$
is known as a subconvex bound. For $m=1$ and 2 such bounds have been established, but for $m\ge3$ the problem remains a central and widely open one.

For the variances $V(\pi)$, upper bounds and even their finiteness is problematic, as is discussed in detail in the paper. For $m$ fixed we expect that $V(\pi)$ is of order $\log c(\pi)$ . For various families of $\pi$'s that are studied in \cite{SST} one can show that the typical member has $V(\pi)$ of order $\log c(\pi)$, but individual growing lower bounds seem difficult to prove, even conditionally assuming the Riemann hypothesis. The main result of the paper is the \textit{universal} lower bound for $V(\pi)$ in $(2)$ \footnote{Assume $L(s,\pi)$ has no pole.     } :
$$
V(\pi) > 2\pi \log 2\ \Big|1-L\Big(\frac{1}{2},\pi\Big)\Big|^2.
\leqno(5)
$$
If $L(\frac{1}{2},\pi)=0$ this recovers the form that is stated and proved in the paper (the proof  is the same and uses (3) after restricting the $x$ integral to $[1,2]$).

Similarly one can prove a universal lower bound for the mean---square .

For any $\pi$ :
$$
I(\pi) =\displaystyle \int_{-\infty}^\infty \frac{|L(\frac{1}{2}+it,\pi)|^2}{\frac{1}{4}+t^2} dt > \pi
\leqno(6)
$$
(the right-hand side is $\frac{\pi}{2}$ if the L-function has a pole.)
This follows from the following Parseval identity which generalizes (12.5.4) in Titchmarsh \cite{Ti} (with the usual meaning if either is infinite)
$$
\mathrm{LHS\ of\ (6)} = 2\pi  \int_0^\infty \Big|\sum_{n\le x} \lambda_\pi(n) -\delta_\pi x \Big|^2 \frac{dx}{x^2}
\leqno(7)
$$
where $\delta_{\pi}=1$ if $m=1$  and $\pi$ is the trivial representation, and $0$ otherwise.

Perhaps the most basic question that presents itself is whether
$$
\lim_{\pi \rg \infty}  I(\pi)= \infty
$$
when $\pi$ tends to infinity in a suitable sense.
 


\vskip6mm

\noindent \textbf{Proof of (5) and its variant}

\vskip6mm
\textit{We first assume that $L(s,\pi)$ has no pole.} Let
$$
\mathrm  \qquad J(x)= \left\{
\begin{array}{ccl}
1 &\mathrm{if} &0<x\le1\\
0 &\mathrm{if} &x>1.
\end{array}
\right.
\leqno(8)
$$
We smooth $J$ as a function on the multiplicative group $\R_{>0}$ as follows : fix $\varphi \in C_c^\infty(\R)$, $\varphi \ge 0$, $\int_{-\infty}^\infty \varphi(X)dX=1$, $\varphi$ even, and $\rm Support (\varphi) \subset [-\frac{1}{2},\frac{1}{2}]$. For $\varepsilon>0$ set $\varphi_\varepsilon(x) = \frac{1}{\varepsilon}\varphi(\frac{x}{\varepsilon})$ and
$$
J_\varepsilon(x) = \int_0^\infty J(xy^{-1})\varphi_\varepsilon(\log y)\frac{dy}{y}=(J \star\psi_{\varepsilon})(x)
\leqno(9)
$$
(multiplicative convolution, $\psi_{\varepsilon}(x)= \varphi_{\varepsilon}(\log x).)$
 
$J_\varepsilon(x)$ is smooth on $(0,\infty)$ and satisfies for $\varepsilon$ sufficiently small :
$$
\left.
(A) \hspace{2cm} \begin{array}{lll}
J_\varepsilon(x) =1 &\mathrm{for} &0<x<1-\varepsilon\\
J_\varepsilon(x) =0 &\mathrm{for} &x>1+\varepsilon\\
0\le J_\varepsilon(x) \leq1 &\mathrm{for} &0<x<\infty\\
\end{array}
\right\}
\leqno(10)
$$
$(B)$\kern1cm  $J_\varepsilon(x) \rg J(x)$ uniformly for $x$ outside any neighbourhood of 1.

\noindent $(C)$\kern1cm For Re$(s)>0$,
$$
\tilde{J}_\varepsilon(s) := \Mm J_{\varepsilon}(s) =  \int_0^\infty J_\varepsilon(x) x^s \frac{dx}{x} = \Mm J(s) \Mm \psi_{\varepsilon}(s)= \frac{\hat{\varphi}(i\varepsilon s)}{s}
\leqno(11)
$$
where 
$$
\hat{\varphi}(\xi) = \int_{-\infty}^\infty \varphi(X) e^{-i\xi X}dX.
$$
In particular $\tilde{J}_\varepsilon(s)$ is analytic in $\C$ except for a simple pole at $s=0$, with residue 1 at that point. Moreover from (11) it follows that $\tilde{J}_\varepsilon(s)$ is rapidly decreasing in $|t|$ for $s=\sigma+it$ and uniformly so for 
$$
\sigma_0 \le \sigma \le \sigma_1.
\leqno(12)
$$
For $x>0$ define $H_\varepsilon(x)$ by
$$
H_\varepsilon(x) = \frac{1}{2\pi i} \int\limits_{\mathrm{Re}(s)=2} \Big[L\Big( s+\frac{1}{2},\pi\Big) -L\Big(\frac{1}{2},\pi\Big)\Big] \tilde{J}_\varepsilon(s) x^sds.
\leqno(13)
$$
This integral converges absolutely in view of standard bounds for $L(s,\pi)$ in vertical strips and the rapid decay of $\tilde{J}_\varepsilon(s)$. Moreover the series (1) converges absolutely for Re$(s)=5/2$ (Jacquet-Shalika \cite{JS}) and so we can integrate the series definition of $L(s+\frac{1}{2},\pi)$ in (13) term by term and use the Mellin inversion and (10) to conclude that for $x>0$,
$$
H_\varepsilon(x) = \sum_{n=1}^\infty \frac{\lambda_\pi(n)}{\sqrt{n}}J_\varepsilon\Big(\frac{n}{x}\Big) -L \Big(\frac{1}{2},\pi\Big) J_\varepsilon\Big(\frac{1}{x}\Big).
\leqno(14)
$$

For $x$ bounded, the sum is finite. By (9A) the support is bouded away from 0.
Thus $H_\varepsilon(x)$ is smooth and supported away from 0.

From (13) and shifting the contour to Re$(s)=-A $, a large negative number, picks up no poles (since the pole of $\tilde{J}_\varepsilon(s)$ at $s=0$ is cancelled by the difference of the $L$--values), we see that $H_\varepsilon(x)$ decays faster that $x^{-A}$ for any $A$ as $x\lgr \infty$. Hence $H_\varepsilon(x)$ is smooth on $(0,\infty)$ and in
$$
L^2(\R_{>0}^{\times},\frac{dx}{x}).
\leqno(15)
$$

As we just did, we shift the contour integral in (13) to Re$(s)=0$, which is again justified by standard polynomial bounds in $|t|$  for $L(s,\pi)$ and the rapid decay of $\tilde{J}_\varepsilon(s)$ in $|t|$. This yields
$$
\begin{array}{ll}
H_\varepsilon(x) &= \dfrac{1}{2\pi}\displaystyle \int_{-\infty}^\infty \Big[L\Big(\dfrac{1}{2}+it,\pi\Big)-L\Big(\dfrac{1}{2},\pi\Big)\Big]\tilde{J}_\varepsilon(it)x^{it}dt\\
&= \dfrac{1}{2\pi}\displaystyle \int_{-\infty}^\infty \Big[L\Big(\dfrac{1}{2}+it,\pi\Big)-L\Big(\dfrac{1}{2},\pi\Big)\Big] \dfrac{\hat{\varphi}(\varepsilon t)}{it} x^{it}dt.
\end{array}
\leqno(16)
$$
Hence the rapidly decreasing smooth $H_\varepsilon(x)$ on $\R_{>0}^\times$ satisfies
$$
\Mm H_\varepsilon(-it) = \frac{L(\frac{1}{2}+it,\pi)-L(\frac{1}{2},\pi)}{it}\  \hat{\varphi}(\varepsilon t) \ (t\in \R).
\leqno(17)
$$
We will write $\Ff H(t)= \Mm H(-it)$ for $H$ a fonction on $\R^{\times}_+$.

Now let $\varepsilon\rg 0$. If $\Ff H(t):= \dfrac{L(\frac{1}{2}+it,\pi)-L(\frac{1}{2},\pi)}{it}$ is not in $L^2(\R)$ there is nothing to prove as (5) is then valid trivially.

So we can assume that
$$
\Ff H(t) \in L^2(\R).
\leqno(18)
$$
Now from (17) and (9) we have that
$$
\Ff H_\varepsilon(t) \lgr \Ff H(t) \ \textrm{uniformly\ on\ compacta}
\leqno(19)
$$
and
$$
H_\varepsilon(x) \lgr H(x), \textrm{uniformly\ on\ compacta\ in\ }(0,\infty)\ba \N,
$$

where
$$
H(x) := \left\{
\begin{array}{cl}
\displaystyle\sum_{n\in x} \dfrac{\lambda\pi(n)}{\sqrt{n}} -L( \frac{1}{2},\pi),\ &x>1\\
0 &\textrm{for}\ x<1.
\end{array}
\right.
\leqno(20)
$$
Moreover from (17) and that $\Ff H\in L^2(\R)$, it follows from the dominated convergence theorem (or more simply by estimating the tails of the $t$ integral uniformly) that $\Ff H_{\varepsilon} \lgr \Ff H$ in $L^2(\R)$. Hence by Parseval $H_\varepsilon$ 
converges in $L^2(\R_{>0}^\times,\frac{dx}{x})$ and from (20) the limit of $H_\varepsilon$ must be $H$. Thus $H\in L^2(\R^\times_{>0},\frac{dx}{x})$ and $\Ff H(t)$ and $H(x)$ are Fourier pairs and hence by Parseval, we have that
$$
V(\frac{1}{2},\pi)=\vert\vert \Ff H \vert\vert_2^2 = 2\pi \int_0^\infty | H(x)|^2 \frac{dx}{x}.
\leqno(21)
$$

This proves (3).Now following the argument in the paper we note that on $[1,2]\subset \R_{>0}^\times$, $H(x) = 1-L( \frac{1}{2},\pi)$ so that the R.H.S. of (21) is at least
$$
2\pi |1-L(\frac{1}{2},\pi)|^2 \int_1^2 \frac{dx}{x} = 2\pi \log2|1-L(\frac{1}{2},\pi)|^2.
\leqno(22)
$$
Moreover for $x>2$ it is easy to see that $H(x)\equiv 0$ is impossible, and hence we arrive at (5).

We note that the argument obviously extends to the case where $\pi= \pi_1 \times \pi_2 \times...\times\pi_r$ as in the Introduction to the main text, and that we have not used the Ramanujan hypothesis. Moreover, by an obvious change of variable (of the form $s\rightarrow s+i\tau$), this implies Theorem A for any $s$ on the critical line.

Now consider Theorem B. We cannot use again the translation argument just introduced : this shifts the pole at $(s=1)$ to a pole at $s=1+i\tau$. However, the computation is the same and we extend the previous proof, assuming the pole is at $(s=1)$.

We define $H_\varepsilon$ as before by (13) ; the  equality (14) remains the same. However, when we shift the integral to $\mathrm{Re}(s)=A\ll 0$, we pick up a pole at $s=1/2$, with residue $\kappa\,\tilde{J}_\varepsilon(1/2)x^{1/2}$. Thus the properties of smoothness and decrease are now true for
$$
H_\varepsilon^1(x) = H_\varepsilon(x) - \kappa\ \tilde{J}_\varepsilon(1/2) x^{1/2}.
\leqno(23)
$$

Furthermore, this pole also occurs when we shift the integral to $\mathrm{Re}(s)=0$, so the equality (16) is now true for $H_\varepsilon^1$. The some arguments now show that
$$ 
\Ff\, H^1 = \frac{L(1/2+it,\pi)-L(1/2,\pi)}{it},
$$
 where
$$
H^1(x) = H(x) - \frac{\kappa}{1/2} x^{1/2}
$$
since $\tilde{J}_\varepsilon(1/2) = \frac{\hat{\varphi}(1/2\varepsilon)}{1/2}\rg 2$ when $\varepsilon\rg 0$. Finally, $\frac{1}{2\pi} V (1/2, \pi)$ is now bounded below by
$$
\int_0^1 |2\kappa x^{1/2}|^2 \frac{dx}{x} + \int_1^2 |1-L(1/2,\pi) - 2\kappa x^{1/2}|^2\frac{dx}{x},
$$
leading to Theorem B if $L(1/2,\pi)=0$.

\eject

\noindent\textbf{Proof of (6) and its variant}

\vskip2mm

We now consider the new function
$$
H_\varepsilon(x) = \frac{1}{2i\pi} \int_{\Res=2} [L(s,\pi)-L(0,\pi)] \tilde{J}_\varepsilon(s)x^sds.
 \leqno(24)
 $$
 
 For $\mathrm{Re}(s)\ge 2$ we obtain, integrating term by term :
$$
 H_\varepsilon(x) = \sum_{n=1}^\infty \lambda_\pi (n) J_\varepsilon\Big(\frac{n}{x}\Big) -L(0,\pi) J_\varepsilon \Big(\frac{1}{x}\Big).
 \leqno(25)
 $$
 
 The function $(L(s,\pi)-L(0,\pi) \tilde{J}_\varepsilon(s)$ is holomorphic. Shifting as before the integral to $\mathrm{Re}(s)=A$, $A\ll 0$, we find again that $H_\varepsilon$ is smooth and of rapid decrease. We have
$$
 K_\varepsilon(x) :=  x^{-1/2}H_\varepsilon(x) = \frac{1}{2\pi} \int_{-\infty}^{+\infty} \frac{L(1/2+it)-L(0)}{1/2+it} \cdot \hat{\varphi}(\varepsilon t-i\frac{\varepsilon}{2})x^{it}dt.
 \leqno(26)
 $$
 However, $\frac{L(0)}{1/2+it}\in L^2(\R)$. Again, there is nothing to prove if $\frac{L(1/2+it)}{1/2+it}\not\in L^2$, and otherwise the first factor in the integrand of (26) is $L^2$. Now the previous arguments show that $K_\varepsilon(x) \rg K(x)$, uniformly on compacta in $(0,\infty) -\N$,
with
$$
K(x) = \left\{ \begin{array}{l}
x^{-1/2}\Big(\sum\limits_{n\le x}\lambda_\pi (n) -L(0,\pi)\Big),\ x\ge 1\\
0 \qquad \mathrm{for}\ x<1.
\end{array}\right.
\leqno(27)
$$
From (26) it follows that
$$
\Mm K_\varepsilon(-it) \lgr \frac{L(1/2+it)-L(0)}{1/2+it}
$$
in $L^2(\R)$. Moreover, 
$$
\Mm K(-it) =\frac{L(1/2+it)-L(0)}{1/2+it}. 
$$
If $\chi$ is the characteristic function of $(x\ge 1)$, $\Mm(x^{-1/2}\chi)(-it)=\frac{1}{1/2+it}$. Thus $x^{-1/2} \sum\limits_{n\le x}\lambda_\pi(n)$ $(x\ge 1)$ and $\frac{L(1/2+it)}{1/2+it}$ are associated, and this implies
$$
\int_{-\infty}^{+\infty} \Big|\frac{L(1/2+it)}{1/2+it}\Big|^2 dt > 2\pi \int_1^2 x^{-2}dx=\pi.
$$
Finally, Theorem D can be proven by combining this computation with the argument given for Theorem~B, introducing the residue at $(s=1)$.



\eject

\vskip6mm

Laurent Clozel

Math\'{e}matiques, B\^{a}timent 307

Universit\'{e} Paris-Saclay

91405 Orsay Cedex

France

\vskip6mm

Peter Sarnak

School of Mathematics

Institute for Advanced Study

Princeton NJ  08540

USA

\end{document}